\documentclass[11pt]{amsart}

\usepackage{amsmath,amsthm,amssymb,amsxtra,amsfonts,amsbsy,mathrsfs}
\usepackage{verbatim}
\usepackage[all]{xy}
\newtheorem{theorem}{Theorem}[section]
\newtheorem{corollary}[theorem]{Corollary}
\newtheorem{lemma}[theorem]{Lemma}
\newtheorem{proposition}[theorem]{Proposition}
\newtheorem{definition}[theorem]{Definition}

\newtheorem{remark}[theorem]{Remark}
\newtheorem{example}[theorem]{Example}

\newtheorem*{introtheorem}{Theorem \ref{thm:hyp}}
\newtheorem*{introcorollary}{Theorem \ref{cor:tamedegen}}

\newcommand{\Q}{\mathbb{Q}}
\newcommand{\Z}{\mathbb{Z}}
\newcommand{\C}{\mathbb{C}}

\newcommand{\til}[1]{\tilde{#1}}
\newcommand{\mf}[1]{\mathfrak{#1}}
\newcommand{\OO}{\mathcal{O}}
\newcommand{\XX}{\varpi_{\mf{X}/\mf{X}}}
\newcommand{\XS}{\varpi_{\mf{X}/S}}

\newcommand{\XBS}{\Omega_{X_\bullet/S}}

\newcommand{\LXS}{L_{\mf{X}/S}}
\newcommand{\pLXS}{{\bigwedge}^{\!\!s} L_{\mf{X}/S}}
\newcommand{\DRB}{\Omega^{\bullet}_{M^0/S}}
\newcommand{\DR}{\Omega_{M^0/S}}
\newcommand{\HH}{\mathcal{H}}
\newcommand{\XPS}{\varpi_{\mf{X}'/S}}
\newcommand{\XSB}{\varpi_{\mf{X}/S}^\bullet}

\DeclareMathOperator{\Spec}{Spec}

\DeclareMathOperator{\Hom}{Hom}
\DeclareMathOperator{\Ext}{Ext}
\DeclareMathOperator{\Tor}{Tor}

\DeclareMathOperator{\Aut}{Aut}

\SelectTips{xy}{12}

\makeatletter
\@namedef{subjclassname@2010}{ \textup{2010} Mathematics Subject 
Classification} \makeatother

\begin{document}


\title[de Rham Theory for Tame Stacks]{de Rham Theory for Tame Stacks and\\ Schemes with Linearly Reductive 
Singularities}
\author{Matthew Satriano}
\address{Department of Mathematics, University of Michigan, 2074 East Hall, Ann Arbor, MI 48109-1043, USA}
\curraddr{}
\email{satriano@umich.edu}

\date{}

\maketitle

\vspace{-1cm}

\begin{abstract}
We prove that the Hodge-de Rham spectral sequence for smooth proper tame Artin stacks in characteristic $p$ (as defined by 
Abramovich, Olsson, and Vistoli) which lift mod $p^2$ degenerates. We push the result to the coarse spaces of such stacks, 
thereby obtaining a degeneracy result for schemes which are \'etale locally the quotient of a smooth scheme by a finite 
linearly reductive group scheme.
\end{abstract}

\vspace{0.5cm}

Given a scheme $X$ smooth and proper over a field $k$, the cohomology of the algebraic de Rham complex $\Omega^\bullet_{X/k}$ 
is an important invariant of $X$, which, when $k=\C$, recovers the singular cohomology of $X(\C)$.  When the Hodge-de Rham 
spectral sequence
\[ 
E_1^{st}=H^t(\Omega^s_{X/k})\Rightarrow H^n(\Omega^\bullet_{X/k})
\]
degenerates, the invariants 
$\dim_k H^n(\Omega^\bullet_{X/k})$ break up into sums of the finer invariants $\dim_k H^t(\Omega^s_{X/k})$.  The 
degeneracy of this spectral sequence for smooth proper schemes in characteristic 0 was first proved via analytic methods.  
It was not until much later that Faltings \cite{faltings} gave a purely algebraic proof by means of $p$-adic Hodge Theory.  
Soon afterwards, 
Deligne and Illusie \cite{di} gave a substantially simpler algebraic proof by showing that the degeneracy of the 
Hodge-de Rham spectral sequence in characteristic 0 is implied by its degeneracy for smooth proper schemes in 
characteristic $p$ that lift mod $p^2$.  Their method therefore extends de Rham Theory to the class of smooth proper 
schemes in positive characteristic which lift.  A version of de Rham 
Theory also exists for certain singular schemes.  Steenbrink showed \cite[Thm 1.12]{Steenbrink} that if $k$ is a field of 
characteristic 0, $M$ a proper $k$-scheme with quotient singularities, and $j:M^0\hookrightarrow M$ its smooth locus, then 
the hypercohomology spectral sequence 
\[
E_1^{st}=H^t(j_*\Omega_{M^0/k}^s)\Rightarrow H^n(j_*\Omega_{M^0/k}^{\bullet})
\]
degenerates and $H^n(j_*\Omega_{M^0/k}^{\bullet})$ agrees with $H^n(M(\C),\C)$ when $k=\C$.  As we explain in this paper, a 
version of this theorem is true in positive characteristic as well: if $k$ has characteristic $p$ and 
$M$ is proper with quotient singularities by groups whose orders are prime to $p$, then the above spectral sequence 
degenerates for $s+t<p$ provided a certain liftability criterion is satisfied (see Theorem \ref{thm:goodquot} for precise 
hypotheses).\\
\\
As a warm-up for the rest of the paper, we begin by showing how Steenbrink's result can be reproved using the theory of 
stacks.  The idea is as 
follows.  Every scheme $M$ as above is the coarse space of a smooth Deligne-Mumford stack $\mf{X}$ whose stacky structure 
is supported at the singular locus of $M$.  We show that the de Rham cohomology $H^n(\Omega^\bullet_{\mf{X}/k})$ of the stack 
agrees with $H^n(j_*\Omega_{M^0/k}^{\bullet})$.  After checking that the method of Deligne-Illusie extends to Deligne-Mumford 
stacks, we recover Steenbrink's result as a consequence of the degeneracy of the Hodge-de Rham spectral sequence for 
$\mf{X}$.\\
\\
The above extends de Rham Theory to the class of schemes with quotient singularities by groups whose orders are prime to 
the characteristic, but in positive characteristic 
this class of schemes contains certain ``gaps'' and it is natural to ask if de Rham Theory can be extended further.  For 
example, in all characteristics except for 2, the affine quadric cone $\Spec k[x,y,z]/(xy-z^2)$ 
can be realized as the quotient of $\mathbb{A}^2$ by $\Z/2\Z$ under the action $x\mapsto -x$, $y\mapsto -y$.  In 
characteristic 2, however, this action is trivial.  If we allow quotients not just by finite groups, but rather finite 
group schemes, then we can realize the cone as $\mathbb{A}^2/\mu_2$ where $\zeta\in\mu_2(T)$ acts as 
$x\mapsto\zeta x$, $y\mapsto\zeta y$.  This is an example of what we call a scheme with 
\emph{linearly\ reductive\ singularities}; that is, a scheme which is \'etale locally the quotient of a smooth scheme 
by a finite flat linearly reductive group scheme.\\
\\
One of the main results of this paper is that de Rham Theory can be extended to the class of schemes with isolated linearly 
reductive singularities.  As with 
Steenbrink's result, we prove this by passing through stacks.  Just as schemes with quotient singularities are coarse spaces 
of smooth Deligne-Mumford stacks whose stacky structure is supported at the singular locus, we show in Theorem 
\ref{prop:cslrs} that schemes with linearly reductive singularities are coarse spaces of smooth tame (Artin) stacks (as 
introduced in \cite{tame}) whose stacky structure is supported at the singular locus.  So, to 
extend de Rham Theory to schemes with isolated linearly reductive singularities, we first show the degeneracy of a type of 
Hodge-de Rham spectral sequence for tame stacks.  We should emphasize that, unlike in the case of 
Deligne-Mumford stacks, there are technical barriers to extending the method of Deligne-Illusie to Artin stacks or even 
tame stacks, first and foremost being that relative Frobenius does not behave well under smooth base change.  It should 
also be noted that it is \emph{a\ priori} not clear what the definition 
of the de Rham complex of a tame stack $\mf{X}$ should be.  One can use the cotangent complex $L_{\mf{X}}$ of the stack 
(see \cite[\S 15]{lmb} and \cite[\S 8]{sheaves}) to define the derived de Rham complex $\bigwedge^{\bullet}L_{\mf{X}}$; 
alternatively, one can use a more naive sheaf of differentials $\varpi_{\mf{X}}^1$ on the lisse-\'etale site 
of $\mf{X}$ whose restriction to each $U_{et}$ is $\Omega^1_{U}$, for every $U$ smooth over $\mf{X}$.  The latter has
the advantage that it is simpler, but it is \emph{not\ coherent}; the cotangent complex, on the other hand, has 
coherent cohomology sheaves.  We take the naive de Rham complex as our definition, but it is by comparing this complex 
with the derived de Rham complex that we prove our main result for tame stacks:
\begin{introcorollary}
Let $\mf{X}$ be a smooth proper tame stack over a perfect field $k$ of characteristic $p$.  If $\mf{X}$ lifts mod $p^2$, 
then the Hodge-de Rham spectral sequence
\[
E_1^{st}=H^t(\varpi^s_{\mf{X}/k})\Rightarrow H^n(\varpi^\bullet_{\mf{X}/k})
\]
degenerates for $s+t<p$ $($see the Notation section below$)$.
\end{introcorollary}
From Theorem \ref{cor:tamedegen} and Theorem \ref{prop:cslrs}, we are able to deduce
\begin{introtheorem}
Let $M$ be a proper $k$-scheme with isolated linearly reductive singularities, where $k$ is a perfect field of characteristic 
$p$.  Let $j:M^0\hookrightarrow M$ be the smooth locus of $M$ and let $\mf{X}$ be as in Theorem \ref{prop:cslrs}.  If 
$\mf{X}$ lifts mod $p^2$, then the hypercohomology spectral sequence
\[
E_1^{st}=H^t(j_*\Omega^s_{M^0/k})\Rightarrow H^n(j_*\Omega^\bullet_{M^0/k})
\]
degenerates for $s+t<p$.
\end{introtheorem}
We should mention that unlike in the case of quotient singularities, the cohomology groups 
$H^n(\varpi^\bullet_{\mf{X}/k})$ and $H^n(j_*\Omega^\bullet_{M^0/k})$ no longer agree, so some care is needed in showing how 
Theorem \ref{thm:hyp} follows from the degeneracy of the Hodge-de Rham spectral sequence of the stack.\\
\\
It is desirable, of course, to remove the stack from the statement Theorem \ref{thm:hyp}.  We show in Theorem \ref{thm:lift} 
that if the dimension of $M$ is at least 4, 
the liftability of $M$ implies the liftability of $\mf{X}$.  In this case, we therefore have a purely scheme-theoretic 
statement of Theorem \ref{thm:hyp}.  We end the paper by proving a type of Kodaira vanishing theorem within this setting.\\
\\
This paper is organized as follows.  In Section 1, we
%Section 1 is a warm-up for the rest of the paper.  We 
begin by reviewing some 
background material and giving an outline of \cite[Thm 2.1]{di} as some of the technical details will be used later.  We then 
consider de Rham Theory for Deligne-Mumford stacks and show how stacks can be used to recast Steenbrink's result.  
%$\cite{Steenbrink}$ that if $k$ is a characteristic 0 field, $M$ a proper $k$-scheme 
%with quotient singularities, and $j:M^0\rightarrow M$ its smooth locus, then the spectral sequence 
%$E_1^{st}=H^t(j_*\Omega_{M^0/k}^s)\Rightarrow H^n(j_*\Omega_{M^0/k}^{\bullet})$ degenerates.  
The purpose of Section 2 
is to find a way around the problem that the method of Deligne and Illusie does not carry over directly to the 
lisse-\'etale site of Artin stacks.  Since relative Frobenius \emph{does} behave well under \'etale base change, 
our solution is to prove a Deligne-Illusie result on the \'etale site of $X_{\bullet}$, where $X\rightarrow \mf{X}$ 
is a smooth cover of a smooth tame stack by a scheme, and $X_\bullet$ is the simplicial scheme obtained by taking fiber 
products over $\mf{X}$.  The key technical point here is showing that \'etale locally on the coarse space of $\mf{X}$, the 
relative Frobenius for $\mf{X}$ lifts mod $p^2$.  
%Here $\mf{X}$ is any smooth Artin stack which lifts mod $p^2$ and satisfies a mild local structure condition.  
In Section 3, 
we prove that the naive de Rham complex and the derived de Rham complex above compute the same cohomology, and show how 
this result implies the 
degeneracy of the Hodge-de Rham spectral sequence for smooth proper tame stacks which lift mod $p^2$.  In Section 4, 
we prove Theorem \ref{thm:hyp}.\\
%The main focus of Section 4 is to prove Theorem \ref{thm:hyp}
%, although we begin by proving two results which may be of independent interest: 
%the first is Theorem \ref{prop:cslrs}; the second is a type of Cartier isomorphism for schemes with linearly reductive singularities.  
%The proof of Theorem \ref{prop:cslrs} involves generalizing the classical Chevalley-Shephard-Todd Theorem 
%\cite[\S 5 Thm 4]{bourbaki} to the case of finite linearly reductive group schemes.  This is done in Appendix A.\\
\\
$\textbf{Acknowledgements.}$  I would like to thank Dustin Cartwright, Ishai Dan-Cohen, and Anton Geraschenko for helpful 
conversations.  Most of all, I am grateful to my advisor, Martin Olsson, both for his guidance and his help in editing this 
paper.\\
\\
$\textbf{Notation.}$  Unless otherwise stated, all Artin stacks are assumed to have finite diagonal.  
If $\mf{X}$ is an Artin stack over a scheme $S$, we let $\mf{X}'$ denote the pullback of $\mf{X}$ by the absolute 
Frobenius $F_S$.  We usually drop the subscript on the relative Frobenius $F_{\mf{X}/S}$, denoting it by $F$.  
Given a morphism $g:\mf{X}_1\rightarrow\mf{X}_2$ of $S$-stacks, we denote by 
$g':\mf{X}'_1\rightarrow\mf{X}'_2$ the induced morphism.\\
%Given a simplicial topos $X_\bullet$, we let $X_\bullet^{+}$ denote the strictly simplicial topos obtained 
%by forgetting the degeneracy maps.\\
Given a morphism $g:\mf{X}_1\rightarrow\mf{X}_2$ of Artin stacks and complex of sheaves $\mathcal{F}^\bullet$ on $\mf{X}_1$, 
we do not use the shorthand $g_*\mathcal{F}^\bullet$ when we mean $Rg_*\mathcal{F}^\bullet$.  For us, 
$g_*\mathcal{F}^\bullet$ always denotes the complex obtained by applying the functor $g_*$ to the complex 
$\mathcal{F}^\bullet$.\\
Lastly, we say a first quadrant spectral sequence $E_{r_0}$ ``degenerates for $s+t<N$'' if for all $r\geq r_0$ and 
all $s$ and $t$ satisfying $s+t<N$, all of the differentials to and from the $E_r^{st}$ are zero.

%If $k$ is a field of positive characteristic, it is standard to denote $\character k$ by $p$; however, 
%it is also standard to use $p$ and $q$ for spectral sequence indices.  We are often interested in proving the degeneracy 
%of a spectral sequence for $p+q<\character k$.  Writing $p+q<p$ would, of course, be confusing, but we hope the phrase 
%``$p+q<\character k$'' causes the reader no confusion.\\
%SHOULD I JUST USE $\ell$ FOR THE CHARACTERISTIC?\\

\section{Steenbrink's Result via Stacks}
\label{sec:viastacks}
\subsection{Review of Deligne-Illusie}
\label{subsec:review}
We briefly review the proof of \cite[Thm 2.1]{di} and explain how it generalizes to Deligne-Mumford stacks.  Having an 
outline of this proof will be useful for us in Section $\ref{sec:simplicial}$.\\
%For further background material, see $\cite{katz}$.\\
\\
Let $S=\Spec k$ be a perfect field of characteristic $p$.  For any $S$-scheme $X$, let $F_X:X\rightarrow X$ 
be the absolute Frobenius, which acts as the identity on topological spaces and sends a local section $s\in\OO_X(U)$ 
to $s^p$.  We have the following commutative diagram, where $F_{X/S}$ is the relative Frobenius and 
the square is cartesian 
\[
\xymatrix{
X\ar[dr]\ar[r]^{F_{X/S}}\ar@/^2.5pc/[rr]^{F_X} & X'\ar[r]\ar[d] & X\ar[d]\\
 & S\ar[r]^{F_S} & S
}
\]
We drop the subscript on the relative Frobenius $F_{X/S}$, denoting it by $F$.  If $X$ is locally of finite type over $S$, so 
that it is locally $\Spec k[x_1,\dots,x_n]/(f_1,\dots,f_m)$, where $f_j=\sum a_{j,I}x^I$, then $X'$ is locally 
$\Spec k[x_1,\dots,x_n]/(f_1^{(p)},\dots,f_m^{(p)})$, where $f_j^{(p)}:=\sum a_{j,I}^p x^I$.  The relative Frobenius 
morphism is then given by sending $x_i$ to $x_i^p$ and $a\in k$ to $a$.\\
\\
Our primary object of study is the de Rham complex $\Omega^\bullet_{X/S}$.  
%The maps $d:\Omega_{X/S}^{k}\rightarrow \Omega_{X/S}^{k+1}$ in the complex are defined as the composite \[\Omega_{X/S}^{k}=\OO_X\otimes_{\OO_X}\Omega_{X/S}^{k}\stackrel{d\otimes id}{\longrightarrow}\Omega_{X/S}^{1}\otimes_{\OO_X}\Omega_{X/S}^{k}\longrightarrow\Omega_{X/S}^{k+1},\]where the last map is given by $\eta\otimes\omega\mapsto\eta\wedge\omega$.  Note that the maps in the complex $\Omega^\bullet_{X/S}$ 
The maps in this complex are not $\OO_X$-linear.  To 
correct this ``problem'' we instead consider $F_*\Omega^\bullet_{X/S}$, whose maps $\emph{are}$ 
$\OO_{X'}$-linear.  It is now reasonable to ask how the cohomology of this new complex compares with the cohomology of the 
de Rham complex on $X'$.  An answer is given by:
\begin{theorem}[Cartier isomorphism]
If $X$ is smooth over $S$, then 
%$\HH^i(F_*\Omega_{X/S}^\bullet)=\Omega^i_{X'/S}$.  Moreover, 
there is a unique isomorphism of $\OO_{X'}$-graded algebras
\[
C^{-1}:\bigoplus\Omega^i_{X'/S}\longrightarrow\bigoplus\HH^i(F_*\Omega^\bullet_{X/S})
\]
such that $C^{-1}d(x\otimes 1)$ is the class of $x^{p-1}dx$ for all local sections x of $\OO_{X'}$.
\end{theorem}
Note that once $C^{-1}$ is shown to exist, uniqueness is automatic.  For a proof of this theorem, see \cite[Thm 7.2]{katz}.\\
\\
We are now ready to discuss \cite[Thm 2.1]{di}.
\begin{theorem}
\label{thm:maindi}
Let $W_2(k)$ be the ring of truncated Witt vectors and let $\til{S}=\Spec W_2(k)$.  
If $X$ is smooth over $S$, then to every smooth lift $\til{X}$ of $X$ to $\til{S}$, there is an 
associated isomorphism 
\[
\varphi:\bigoplus_{i<p}\Omega^i_{X'/S}[-i]\longrightarrow \tau_{<p}F_*\Omega_{X/S}^\bullet
\]
in the derived category of $\OO_{X'}$-modules such that $\HH^i(\varphi)=C^{-1}$ for all $i<p$.
\end{theorem}
We give a sketch of the argument.  To define $\varphi$, we need only define 
$\varphi^i:\Omega^i_{X'/S}[-i]\rightarrow \tau_{<p}F_*\Omega_{X/S}^\bullet$ such that $\HH^i(\varphi)=C^{-1}$ 
for each $i<p$.  We take $\varphi^0$ to be the composite 
\[
\OO_{X'}\stackrel{C^{-1}}{\longrightarrow} \HH^0F_*\Omega_{X/S}^\bullet \longrightarrow F_*\Omega_{X/S}^\bullet.
\]
Suppose for the moment that $\varphi^1$ has already been defined.  For $i>1$, we can then define $\varphi^i$ to be 
the composite 
\[
\Omega^i_{X/S}[-i]\stackrel{a[-i]}{\longrightarrow} 
(\Omega^1_{X/S})^{\otimes i}[-i]\stackrel{(\varphi^1)^{\otimes i}}{\longrightarrow} (F_*\Omega^\bullet_{X/S})^{\otimes i}
\stackrel{b}{\longrightarrow} F_*\Omega_{X/S}^\bullet,
\]
where
\[
a(\omega_1\wedge\dots\wedge\omega_i)=\frac{1}{i!}\sum_{\sigma\in S_i}(\textrm{sign}\,\sigma)\,
\omega_{\sigma(1)}\otimes\dots\otimes\omega_{\sigma(i)}
\]
and $b(\omega_1\otimes\dots\otimes\omega_i)=\omega_1\wedge\dots\wedge\omega_i$.\\
\\
Thus, we are reduced to defining $\varphi^1$.  Suppose first that Frobenius lifts; that is, there exists 
$\til{F}$ filling in the diagram 
\[
\xymatrix{
X\ar[r]\ar[d]^{F}\ar@/_1pc/[dd] & \til{X}\ar@{-->}_{\til{F}}[d]\ar@/^1pc/[dd]\\
X'\ar[r]\ar[d] & \til{X}'\ar[d]\\
S\ar[r] & \til{S}
}
\]
where $\til{X}'=\til{X}\times_{\til{S},\sigma} \til{S}$ and $\sigma$ is the Witt vector Frobenius automorphism.  
Let $\mathbf{p}:\OO_X\stackrel{\simeq}{\rightarrow} p\OO_{\til{X}}$ be the morphism sending $x_0$ to $px$ for 
any local section $x$ of $\OO_{\til{X}}$ reducing mod $p$ to $x_0$.  Note that if $x\otimes 1$ is a local section 
of $\OO_{\til{X}}\otimes_{W_2(k),\sigma} W_2(k)=\OO_{\til{X}'}$, then 
$\til{F}^*(x\otimes 1)=x^p+\mathbf{p}(u(x))$ for a unique local section $u(x)$ of $\OO_X$.  We define a 
morphism $f:\Omega^1_{X'/S}\rightarrow F_*\Omega^1_{X/S}$ by 
\[f(dx_0\otimes 1)=x_0^{p-1}dx_0+du(x).
\]
Deligne and 
Illusie show that $\varphi^1$ can be taken to be $f$.  Given two different choices $\til{F_1}$ and $\til{F_2}$ of 
$F$, we obtain a homotopy $h_{12}$ relating $f_1$ and $f_2$, defined by $h_{12}(dx_0\otimes1)=u_2(x)-u_1(x)$.\\
\\
Note that $F$ lifts locally since the obstruction to lifting it lies in
\[
\Ext^1(F^*\Omega^1_{X'/S},\OO_X)=H^1(X,F^*T_{X'/S}).
\]
So, to define $\varphi^1$ in general, we need only patch together the local choices.  This is done as follows.  Let 
$\mathcal{U}=\{U_i\}$ be a cover on which Frobenius lifts and let $\check{\mathcal{C}}^\bullet(\mathcal{U},\mathcal{F})$ 
denote the sheafified version of the $\check{\textrm{C}}$ech complex of a sheaf $\mathcal{F}$.  We define $\varphi^1$ to be 
the morphism in the derived category 
\[
\Omega^1_{X'/S}[-1]\stackrel{\Phi}{\longrightarrow}
\textrm{Tot}(F_*\check{\mathcal{C}}^\bullet(\mathcal{U},\Omega^\bullet_{X/S}))\stackrel{\simeq}{\longleftarrow}
F_*\Omega^\bullet_{X/S},
\]
where
\[
\Phi=(\Phi_1,\Phi_2):\Omega^1_{X'/S}\rightarrow F_*\check{\mathcal{C}}^1(\mathcal{U},\OO_X)
\oplus F_*\check{\mathcal{C}}^0(\mathcal{U},\Omega^\bullet_{X/S})
\]
is given by $(\Phi_1(\omega))_{ij}=h_{ij}(\omega|U'_{ij})$ and $(\Phi_2(\omega))_{i}=f_{i}(\omega|U'_{i})$.  Deligne and 
Illusie further show that this is independent of the choice of covering.  This completes the proof.
\begin{remark}
\label{rmk:morphcomplex}
\emph{In the local case where Frobenius lifts, $\varphi$ is a morphism of complexes.  It is only when patching together 
the local choices that we need to pass to the derived category.}
\end{remark}
\begin{remark}
\label{rmk:didm}
\emph{Using the fact that for any \'etale morphism $g:Y\rightarrow Z$ of $S$-schemes, the pullback of 
$F_{Z/S}:Z\rightarrow Z'$ by $g$ is $F_{Y/S}$, 
%\[
%\xymatrix{
%Y\ar[r]^{F}\ar[d]_{g} & Y'\ar[d]^{g'}\\
%Z\ar[r]^{F} & Z'
%}
%\]
one can check that the proof of Theorem $\ref{thm:maindi}$ works when $X$ is a Deligne-Mumford stack.  
Alternatively, this will follow from the proof of Theorem \ref{thm:disimp} below.}
\end{remark}
Given any abelian category $\mathcal{A}$ with enough injectives, a left exact functor 
$G:\mathcal{A}\rightarrow\mathcal{B}$ to another abelian category, and a bounded below complex of objects $A^\bullet$ of 
$\mathcal{A}$, we obtain a hypercohomology spectral sequence 
\[
E_1^{st}=R^tG(A^s)\Rightarrow R^nG(A^\bullet).
\]
If $\mf{X}$ is a Deligne-Mumford stack over a scheme $Y$, the hypercohomology spectral sequence 
$E_1^{st}=H^t(\Omega^s_{\mf{X}/Y})\Rightarrow H^n(\Omega^\bullet_{\mf{X}/Y})$ obtained in this way is called the 
Hodge-de Rham spectral sequence.\\
\\
As Deligne and Illusie show, Theorem \ref{thm:maindi} implies the degeneracy 
of the Hodge-de Rham spectral sequence for smooth proper schemes.  We reproduce their proof, which requires no 
modification to handle the case of Deligne-Mumford stacks, after first isolating the 
following useful fact from homological algebra.
\begin{lemma}
\label{l:homological}
Let $K$ be a field and $r_0$ a positive integer.  Let $E_{r_0}^{st}\Rightarrow E^{s+t}$ be a first quadrant spectral 
sequence whose terms are finite-dimensional $K$-vector spaces and whose morphisms are $K$-linear.  If $n$ is a non-negative 
integer and 
\[
\sum_{s+t=n} \dim_K E_{r_0}^{st}=\dim_K E^n,
\]
then for all $r\geq r_0$ the differentials to and from the $E_r^{s,n-s}$ are zero.  Hence, if the above equality holds for 
all $n<N$, then the spectral sequence degenerates for $s+t<N$.
\end{lemma}
\begin{proof}
Note that for all $r\geq r_0$
\[
\sum_{s+t=n} \dim_K E_{r+1}^{st}\leq \sum_{s+t=n} \dim_K E_{r}^{st}
\]
with equality if and only if all of the differentials to and from the $E_r^{s,n-s}$ are zero.  Hence 
\[
\sum_{s+t=n} \dim_K E_{\infty}^{st}\leq \sum_{s+t=n} \dim_K E_{r_0}^{st}
\]
with equality if and only if the differentials to and from the $E_r^{s,n-s}$ are zero for all $r\geq r_0$.  Since the 
$E_\infty$ terms are $K$-vector spaces, the extension problem is trivial, and so 
\[
\dim_K E^n=\sum_{s+t=n} \dim_K E_{\infty}^{st}\leq \sum_{s+t=n} \dim_K E_{r_0}^{st}=\dim_K E^n,
\]
which completes the proof.
\end{proof}
\begin{corollary}[{\cite[Cor 2.5]{di}}]
\label{cor:dmdegen}
If $\mf{X}$ is a Deligne-Mumford stack over $S$, which is smooth, proper, and lifts mod $p^2$, then the 
Hodge-de Rham spectral sequence 
\[
E_1^{st}=H^t(\Omega^s_{\mf{X}/S})\Rightarrow H^n(\Omega^\bullet_{\mf{X}/S})
\]
degenerates for $s+t<p$.
\end{corollary}
\begin{proof}
By Theorem \ref{thm:maindi} and Remark \ref{rmk:didm}, we have an isomorphism
\[
\bigoplus_{s<p}\Omega^s_{\mf{X}'/S}[-s]\longrightarrow \tau_{<p}F_*\Omega_{\mf{X}/S}^\bullet
\]
in the derived category of $\OO_{\mf{X}'}$-modules.  It follows that for all $n<p$,
\[
\bigoplus_{s+t=n} H^t(\Omega^s_{\mf{X}'/S})=H^n(\Omega_{\mf{X}/S}^\bullet).
\]
Using the fact that $H^t(\Omega^s_{\mf{X}'/S})=H^t(\Omega^s_{\mf{X}/S})\otimes_{k,F_k}k$, we see
\[
\sum_{s+t=n}\dim_k H^t(\Omega^s_{\mf{X}/S})=\sum_{s+t=n}\dim_k H^t(\Omega^s_{\mf{X}'/S})=
\dim_k H^n(\Omega_{\mf{X}/S}^\bullet),
\]
which, by Lemma \ref{l:homological}, proves the degeneracy of the spectral sequence.
\end{proof}
Deligne and Illusie further show that the degeneracy of the Hodge-de Rham spectral sequence in positive characteristic 
implies the degeneracy in characteristic 0.  While its degeneration in characteristic 0 had previously been known by 
analytic means, this provided a purely algebraic proof.
\begin{corollary}
\label{cor:dmzero}
%[{\cite[Cor 2.7]{di}}]
Let $\mf{X}$ be a Deligne-Mumford stack which is smooth and proper over a field $K$ of characteristic $0$.  Then 
the Hodge-de Rham spectral sequence 
\[
E_1^{st}=H^t(\Omega^s_{\mf{X}/K})\Rightarrow H^n(\Omega^\bullet_{\mf{X}/K})
\]
degenerates.
\end{corollary}
The proof given in \cite[Cor 2.7]{di} for schemes requires only a minor modification.  It uses that 
if $X$ is a smooth proper scheme over a field $K$ of characteristic $0$, then there is an integral 
domain $A$ of finite type over $\Z$, a morphism $A\rightarrow K$, and a smooth proper scheme $Y$ over 
$\Spec A$ which pulls back over $\Spec K$ to $X$.  Since this statement remains true when we allow 
$X$ and $Y$ to be Deligne-Mumford stacks (\cite[p.2]{kv}), the proof given in \cite[Cor 2.7]{di} implies 
Corollary \ref{cor:dmzero} above.

\subsection{de Rham Theory for Schemes with Quotient Singularities}
Let $k$ be a field of characteristic 0 and let $S=\Spec k$.
\begin{definition}
\emph{An $S$-scheme $M$ (necessarily normal) is said to have $\emph{quotient\ singularities}$ if there is an \'etale cover 
$\{U_i/G_i\rightarrow M\}$, where the $U_i$ are smooth over $S$ and the $G_i$ are finite groups.}
\end{definition}
Our goal in this subsection is to use stacks to reprove \cite[Thm 1.12]{Steenbrink} which states:
\begin{theorem}
\label{thm:st}
Let $M$ be a proper $S$-scheme with quotient singularities, and let $j:M^0\rightarrow M$ be its smooth locus.  Then 
the hypercohomology spectral sequence 
\[
E_1^{st}=H^t(j_*\Omega_{M^0/S}^s)\Rightarrow H^n(j_*\Omega_{M^0/S}^\bullet)
\]
of the complex $j_*\Omega^\bullet_{M^0/S}$ degenerates.  Furthermore, if $k=\C$, then 
$H^n(j_*\Omega_{M^0/S}^\bullet)$ agrees with the Betti cohomology $H^n(M(\C),\C)$ of $M$.
\end{theorem}
The following proposition gives the relationship between Deligne-Mumford stacks and schemes with quotient singularities.
\begin{proposition}
\label{prop:dmquot}
Let $M$ be an $S$-scheme and let $j:M^0\rightarrow M$ be its smooth locus.  Then $M$ has quotient singularities 
if and only if it is the coarse space of a smooth Deligne-Mumford stack $\mf{X}$ such that $f^0$ in the diagram
\[
\xymatrix{
\mf{X}^0\ar[r]^{j^0}\ar[d]_{f^0} & \mf{X}\ar[d]^f \\
M^0\ar[r]_j & M
}
\]
is an isomorphism, where $\mf{X}^0=M^0\times_{M}\mf{X}$.
%Then we have a diagram
%\[
%\xymatrix{
% & \mf{X}\ar[d]^f \\
%M^0\ar[r]_j\ar[ur]^{\til{\jmath}} & M
%}
%\]
%where $\til{\jmath}$ is defined as $j^0\circ(f^0)^{-1}$ from the diagram above.
\end{proposition}
For a proof, see \cite[Rmk 4.9]{fant} or \cite[Prop 2.8]{int}.  Vistoli's proposition is slightly 
more general than the proposition above.\\
\\
We give the proof of Theorem $\ref{thm:st}$ after first proving a lemma which compares $j_*\Omega_{M^0/S}^\bullet$ to the 
de Rham complex of a Deligne-Mumford stack.
\begin{lemma}
\label{l:diffqs}
If $M$ is an $S$-scheme with quotient singularities and $\mf{X}$ is as in Proposition \ref{prop:dmquot}, then 
\[
j_*\Omega_{M^0/S}^\bullet=f_*\Omega_{\mf{X}/S}^\bullet.
\]
\end{lemma}
\begin{proof}
To prove this equality, we need only show 
%$\til{\jmath}_*\Omega_{M^0/S}^\bullet=\Omega_{\mf{X}/S}^\bullet$; that is, we want to show 
$j^0_*\Omega_{\mf{X}^0/S}^\bullet=\Omega_{\mf{X}/S}^\bullet$.  That is, given an \'etale morphism 
$U\rightarrow\mf{X}$, we want 
to show $i_*\Omega_{U^0/S}^\bullet=\Omega_{U/S}^\bullet$, where $U^0:=M^0\times_M U$ and $i$ is the projection to $U$.  
Since $\Omega^k_{U/S}$ is locally free, hence reflexive, the following lemma completes the proof.
\end{proof}
\begin{lemma}
\label{l:normal}
Let X be a normal scheme and $i:U\hookrightarrow X$ an open subscheme whose complement has codimension at least 2.  
If $\mathcal{F}$ is a reflexive sheaf on $X$, then the adjunction map $\mathcal{F}\rightarrow i_*i^*\mathcal{F}$ 
is an isomorphism.
\end{lemma}
\begin{proof}
Since $\mathcal{F}$ is reflexive, $\mathcal{F}=\mathcal{H}om(\mathcal{G},\OO_X)$, where 
$\mathcal{G}=\mathcal{F}^{\vee}$.  Therefore,
\[
i_*i^*\mathcal{F}=i_*\mathcal{H}om(i^*\mathcal{G},\OO_U)=\mathcal{H}om(\mathcal{G},i_*\OO_U)
\]
and since $X$ is normal, $i_*\OO_U=\OO_X$.
\end{proof}
\begin{proof}[Proof of Theorem \ref{thm:st}]
Let $\mf{X}$ be as in Proposition \ref{prop:dmquot}.  From Lemma \ref{l:diffqs}, we see that
$j_*\Omega_{M^0/S}^\bullet=f_*\Omega_{\mf{X}/S}^\bullet$ and $j_*\Omega_{M^0/S}^s=f_*\Omega_{\mf{X}/S}^s$ for all $s$.  Since 
the $\Omega_{\mf{X}/S}^s$ are coherent, it follows from \cite[Lemma 2.3.4]{compact} that 
\[
j_*\Omega_{M^0/S}^\bullet=Rf_*\Omega_{\mf{X}/S}^\bullet
\quad 
\textrm{and}
\quad
j_*\Omega_{M^0/S}^s=Rf_*\Omega_{\mf{X}/S}^s.
\]
We see then that
\[
H^n(j_*\Omega_{M^0/S}^\bullet)=H^n(\Omega_{\mf{X}/S}^\bullet)
\quad 
\textrm{and}
\quad
H^t(j_*\Omega_{M^0/S}^s)=H^t(\Omega_{\mf{X}/S}^s).
\]
The proof of Proposition \ref{prop:dmquot} given in \cite[Prop 2.8]{int} shows that $f$ is proper, and so the Hodge-de Rham 
spectral sequence for $\mf{X}$ degenerates by Corollary \ref{cor:dmzero}.  It follows that 
\[
\sum_{s+t=n}\dim_k H^t(j_*\Omega_{M^0/S}^s)=\sum_{s+t=n}\dim_k H^t(\Omega_{\mf{X}/S}^s)=
\dim_k H^n(\Omega_{\mf{X}/S}^\bullet)=\dim_k H^n(j_*\Omega_{M^0/S}^\bullet),
\]
which, by Lemma \ref{l:homological}, proves the degeneracy of the hypercohomology spectral sequence for 
$j_*\Omega_{M^0/S}^\bullet$.\\
\\
We now show that if $k=\C$, then $H^n(j_*\Omega_{M^0/S}^\bullet)=H^n(M(\C),\C)$.  We have shown 
$H^n(j_*\Omega_{M^0/S}^\bullet)=H^n(\Omega^\bullet_{\mf{X}})$, and GAGA for Deligne-Mumford stacks (\cite[Thm 5.10]{gaga}) 
shows
\[
H^n(\Omega^\bullet_{\mf{X}})=H^n(\Omega^\bullet_{\mf{X}^{an}}),
\]
where $\mf{X}^{an}$ is defined in \cite[Def 5.6]{gaga}.  
Note that $\C\rightarrow\Omega^\bullet_{\mf{X}^{an}}$ is a quasi-isomorphism since this can be checked \'etale locally.  It 
follows that 
\[
H^n(\Omega^\bullet_{\mf{X}^{an}})=H^n(\mf{X}^{an},\C).
\]
Lastly, the singular cohomology of $\mf{X}^{an}$ and that of its coarse space, $M(\C)$, are the same.  This is shown in 
\cite[Prop 36]{behrend} for topological Deligne-Mumford stacks with $\Q$-coefficients, but the proof works equally well 
in our situation once it is combined with \cite[Prop 5.7]{gaga}, which states $[U^{an}/G]=[U/G]^{an}$.
\end{proof}
We end this section with some remarks about the situation in positive characteristic.  Suppose $k$ is a perfect field of 
characteristic $p$ and let $S=\Spec k$.
\begin{definition}
\emph{We say an $S$-scheme $M$ (necessarily normal) has $\emph{good}$ $\emph{quotient}$ $\emph{singularities}$ if it has 
an \'etale cover 
$\{U_i/G_i\rightarrow M\}$, where the $U_i$ are smooth over $S$ and the $G_i$ are finite groups of order prime to $p$.}
\end{definition}
Both the proof in \cite{fant} and in \cite{int} (along with Vistoli's Remark 2.9) cited above work in positive 
characteristic.  So, we have the following generalization of Proposition \ref{prop:dmquot}.
\begin{proposition}
\label{prop:dmgoodquot}
Let $M$ be an $S$-scheme, and let $j:M^0\rightarrow M$ be its smooth locus.  Then $M$ has good quotient singularities 
if and only if it is the coarse space of a smooth tame Deligne-Mumford stack $\mf{X}$ $($\cite[\emph{Def\ 2.3.1}]{compact}$)$ 
such that $f^0$ in the diagram 
\[
\xymatrix{
\mf{X}^0\ar[r]^{j^0}\ar[d]_{f^0} & \mf{X}\ar[d]^f \\
M^0\ar[r]_j & M
}
\]
is an isomorphism, where $\mf{X}^0=M^0\times_{M}\mf{X}$.
\end{proposition}
If $\mf{X}$ is a smooth proper tame Deligne-Mumford stack, then the Hodge-de Rham spectral sequence for $\mf{X}$ degenerates 
by Corollary \ref{cor:dmdegen}, and $f_*\mathcal{F}=Rf_*\mathcal{F}$ for any quasi-coherent sheaf on $\mf{X}$ by 
\cite[Lemma 2.3.4]{compact}.  The proof of Theorem \ref{thm:st} therefore gives the following result as well.
\begin{theorem}
\label{thm:goodquot}
Let $M$ be a proper $S$-scheme with good quotient singularities, and let $j:M^0\rightarrow M$ be its smooth locus.  
If $\mf{X}$, as in Proposition \ref{prop:dmgoodquot}, lifts mod $p^2$, then
\[
E_1^{st}=H^t(j_*\Omega_{M^0/S}^s)\Rightarrow H^n(j_*\Omega_{M^0/S}^\bullet)
\]
degenerates for $s+t<p$.
\end{theorem}
As will follow from Theorem $\ref{thm:lift}$ below, if $M$ has dimension at least 4, lifts mod $p^2$, and has isolated 
singularities, then $\mf{X}$ automatically lifts mod $p^2$.

\section{Deligne-Illusie for Simplicial Schemes}
\label{sec:simplicial}
Let $k$ be a perfect field of characteristic $p$ and let $S=\Spec k$.  In this section, we prove a Deligne-Illusie 
result at the simplicial level.  To do so, we must first make sense of the Cartier isomorphism for simplicial schemes.
%Then we have the following commutative diagram
%\[
%\xymatrix{
%Y \ar[rr]^{F_{Y/S}} \ar[dr]_{\rho} & & Y' \ar[r]\ar[d]^{\rho'} & Y\ar[d]^{\rho}\\
% & X\ar[dr]\ar[r]^{F_{X/S}} & X'\ar[r]\ar[d] & X\ar[d]\\
% & & S\ar[r]^{F_S} & S
%}
%\]
%where the squares are cartesian and the composition of the horizontal morphisms are absolute Frobenii.  
%We drop subscripts on the relative Frobenii.
%and often write $\rho$ for both $\rho$ and $\rho'$.
\begin{lemma}
\label{l:Cartier}
Let $X$ and $Y$ be smooth schemes over $S$ and let $\rho:X\rightarrow Y$ be a morphism of $S$-schemes.  
If $C^{-1}$ denotes the Cartier isomorphism, then the following diagram commutes
\[
\xymatrix{
\rho'^*\Omega^i_{X'/S}\ar[r]^-{\rho'^*C^{-1}}\ar[d] & \rho'^*\mathcal{H}^i(F_*\Omega^\bullet_{X/S})\ar[d]\\
\Omega^i_{Y'/S}\ar[r]^-{C^{-1}} & \mathcal{H}^i(F_*\Omega^\bullet_{Y/S})
}
\]
\end{lemma}
\begin{proof}
Using the canonical morphism $\rho'^*\mathcal{H}^i(F_*\Omega^\bullet_{X/S})\rightarrow
\mathcal{H}^i(\rho'^*F_*\Omega^\bullet_{X/S})$ and the multiplicativity property of the Cartier isomorphism, 
we need only check that the diagram commutes for $i=0,1$.  For 
$i=0$, the Cartier isomorphism is simply the kernel map, so the $i=0$ case follows from the 
commutativity of
\[
\xymatrix{
\rho'^*F_*\mathcal{O}_X\ar[r]^-{\rho'^*d}\ar[d] & \rho'^*F_*\Omega^1_{Y/S}\ar[d]\\
F_*\mathcal{O}_Y\ar[r]^-{d} & F_*\Omega^1_{X/S}
}
\]
To handle the $i=1$ case, let $f$ be a local section of $\mathcal{O}_X$ and note that
\[
\xymatrix{
df \ar@{|->}[r] \ar@{|->}[d] & f^{p-1}df\ar@{|->}[d]\\
d(\rho(f))\ar@{|->}[r] & \rho(f)^{p-1}d(\rho(f))
}
\]
\end{proof}
\begin{corollary}
\label{cor:simpcart}
Let $\mathfrak{X}$ be a smooth Artin stack over $S$ 
%with quasi-compact $($not necessarily finite$)$ diagonal, 
and let $X_0\rightarrow\mathfrak{X}$ be a smooth cover by a scheme.  If $X_\bullet$ is the simplicial 
scheme obtained by taking fiber products of $X_0$ over $\mf{X}$, and 
$X'_\bullet$ is its pullback by $F_S$, then there exists a unique isomorphism
\[
C^{-1}: \Omega^i_{X'_\bullet/S}\rightarrow \mathcal{H}^i(F_*\Omega^\bullet_{X_\bullet/S})
\]
such that $C^{-1}(1)=1$, 
$C^{-1}(\omega\wedge\tau)=C^{-1}(\omega)\wedge C^{-1}(\tau)$, and $C^{-1}(df)$ is the 
class of $f^{p-1}df$.
\end{corollary}
\begin{proof}
If such a $C^{-1}$ exists, then its restriction to the $n^{th}$ level of the simplicial scheme is the Cartier 
isomorphism for $X_n$.  Therefore, we need only show existence, which follows from Lemma \ref{l:Cartier}.
\end{proof}
We have now proved the Cartier isomorphism for simplicial schemes.  The other main ingredient in extending Deligne-Illusie 
to simplicial schemes $X_\bullet$, is showing that relative Frobenius for $X_\bullet$ lifts locally.  We note that there are, 
in fact, simplicial schemes for which relative Frobenius does not lift.
\begin{example}
\emph{Let $X_\bullet$ be obtained by taking 
fiber products of $S$ over $B\mathbb{G}_a$.  Lifting Frobenius for $X_\bullet$ is then equivalent to lifting 
Frobenius $F$ of $\mathbb{G}_a$ to a morphism $\til{F}$ of group schemes}
\[
\Spec W_2(k)[x]=\mathbb{G}_{a,\til{S}}\rightarrow\mathbb{G}_{a,\til{S}}=\Spec W_2(k)[x].
\]
\emph{
Since $\til{F}$ reduces to $F$, we must have $\til{F}(x)=x^p+pf(x)$ for some $f(x)\in W_2(k)[x]$.  The condition that 
$\til{F}$ be a group scheme homomorphism implies}
\[
(x+y)^p+pf(x+y)=x^p+y^p+p(f(x)+f(y)), 
\]
\emph{and an easy check shows that this is not possible.}
\end{example}
Although the above example shows that relative Frobenius need not lift locally for an arbitrary simplicial 
scheme, we show that relative Frobenius does lift locally for those simplicial schemes which come from 
smooth tame stacks.  This is the key technical point 
of this section.
\begin{proposition}
\label{l:Artinfroblift}
Let $\mf{X}$ be a smooth tame stack over $S$ with coarse space $M$.  Then \'etale locally on $M$, both $\mf{X}$ and the 
relative Frobenius $F_{\mf{X}/S}$ lift mod $p^2$.
%Let $\mf{X}$ be an Artin stack over $S$ with quasi-compact diagonal.  Further assume that $\mf{X}$ has 
%a smooth affine cover by a smooth affine scheme.  If $\mf{X}$ lifts 
%mod $p^2$, then the relative Frobenius $F_{\mf{X}/S}$ lifts.
\end{proposition}
\begin{proof}
Since the statement of the proposition is \'etale local, by 
%Lemma \ref{l:semidirect} 
\cite[Prop 5.2]{cstlr} we can assume that $M$ is affine and 
$\mf{X}=[U/G]$, where $G=\mathbb{G}_{m,S}^r\rtimes H$ and $H$ is a finite \'etale constant group scheme.  Note that $U$ is 
affine and that the smoothness of $G$ and $\mf{X}$ imply that $U$ is smooth over $S$.\\
\\
As a first step in showing that $\mf{X}$ and $F_{\mf{X}/S}$ lift mod $p^2$, we begin by showing that $BG$ and its relative 
Frobenius lift.  Since the underlying scheme of $G$ is 
$\mathbb{G}_{m,S}^r\times_S H$ and its group structure is determined by the action
\[
H\rightarrow\Aut(\mathbb{G}_m^r)=\Aut(\Z^r),
\]
we can use this same action to define a group scheme $\til{G}=\mathbb{G}_{m,\til{S}}^r\rtimes H$ which lifts $G$.  
It follows that $B\til{G}$ is a lift of $BG$.  Lifting the relative Frobenius of $BG$ is the same as lifting 
the relative Frobenius $F_{G/S}:\mathbb{G}_{m,S}^r\rtimes H\longrightarrow\mathbb{G}_{m,S}^r\rtimes H$ to a group scheme 
homomorphism.  Note that $F_{G/S}$ is given by the identity on $H$ and component-wise multiplication by $p$ on 
$\mathbb{G}_{m,S}^r$.  It therefore has a natural lift mod $p^2$ to the group scheme 
homomorphism given by the identity on $H$ and component-wise multiplication by $p$ on 
$\mathbb{G}_{m,\til{S}}^r$.\\
\\
We now prove that $\mf{X}$ and $F_{\mf{X}/S}$ lift.  There is a natural map $\pi:\mf{X}\rightarrow BG$ which makes
\[
\xymatrix{
U\ar[r]\ar[d] & S\ar[d]\\
\mf{X}\ar[r]^{\pi} & BG
}
\]
a cartesian diagram.  To lift $\mf{X}$ mod $p^2$, it suffices to show that there a stack $\til{\mf{X}}$ and a cartesian 
diagram 
\[
\xymatrix{
\mf{X}\ar[r]\ar[d]_{\pi} & \til{\mf{X}}\ar[d]^{\til{\pi}}\\
BG\ar[r] & B\til{G}
}
\]
The obstruction to the existence of such a diagram lies in $\Ext^2(L_{\mf{X}/BG},\OO_{\mf{X}})$; here $L_{\mf{X}/BG}$ 
denotes the cotangent complex.  Since $\pi$ is representable and smooth, $L_{\mf{X}/BG}$ is a locally free sheaf.
It follows that
\[
R\HH om(L_{\mf{X}/BG},\OO_\mf{X})=\HH om(L_{\mf{X}/BG},\OO_\mf{X}),
\]
which is a quasi-coherent sheaf.  Since $\pi$ is affine and $G$ is linearly reductive, for 
any quasi-coherent sheaf $\mathcal{F}$ on $\mf{X}$, we have 
\[
R\Gamma({\mf{X}},\mathcal{F})=R\Gamma({BG},R\pi_*\mathcal{F})=\Gamma({BG},\pi_*\mathcal{F}).
\]
In particular, 
\[
R\Hom(L_{\mf{X}/BG},\OO_\mf{X})=\Gamma({\mf{X}},\HH om(L_{\mf{X}/BG},\OO_\mf{X}))
\]
and so $\Ext^2(L_{\mf{X}/BG},\OO_\mf{X})=0$.\\
\\
To show that $F_{\mf{X}/S}$ lifts mod $p^2$, it suffices to show that it lifts over our choice $\til{F}_{BG/S}$.  That is, 
it suffices to show that there exists a dotted arrow making the diagram 
\[
\xymatrix{
\mf{X}\ar[d]_{F_{\mf{X}/S}}\ar[r] & \til{\mf{X}}\ar@{-->}[d]\ar@/^1.5pc/[dd]^{\til{F}_{BG/S}\circ\til{\pi}}\\
\mf{X}'\ar[d]_{\pi'}\ar[r] & \til{\mf{X}}'\ar[d]_{\til{\pi}'}\\
BG' \ar[r] & B\til{G}'
}
\]
commute.  The obstruction to finding such a dotted arrow lies in 
$\Ext^1(L_{\mf{X}'/BG'},(F_{\mf{X}/S})_*\OO_\mf{X})$.  As before, we have 
\[
R\HH om(L_{\mf{X}'/BG'},(F_{\mf{X}/S})_*\OO_\mf{X})=\HH om(L_{\mf{X}'/BG'},(F_{\mf{X}/S})_*\OO_\mf{X}),
\]
which is again a quasi-coherent sheaf.  An argument similar to the one above then shows  
$\Ext^1(L_{\mf{X}'/BG'},(F_{\mf{X}/S})_*\OO_\mf{X})=0$, thereby completing the proof.
\end{proof}
We now prove Deligne-Illusie for simplicial schemes coming from smooth tame stacks.
\begin{theorem}
\label{thm:disimp}
Let $\mf{X}$ be a smooth tame stack over $S$.  Let $X_0\rightarrow\mf{X}$ be a smooth cover by a scheme and let $X_\bullet$ 
be the simplicial scheme obtained by taking fiber products of $X_0$ over $\mf{X}$.  
%We further assume that there is a representable \'etale surjection $\mathcal{Y}\rightarrow\mf{X}$ from an 
%Artin stack $\mathcal{Y}$, which has a smooth affine cover by a smooth affine scheme.  
Then, to every lift $\til{X_0}\rightarrow\tilde{\mathfrak{X}}$ of $X_0\rightarrow\mathfrak{X}$, there is a canonically 
associated isomorphism
\[
\varphi: \bigoplus_{i<p}\Omega^i_{X'_\bullet/S}[-i]\rightarrow
\tau_{<p}F_*\Omega^\bullet_{X_\bullet/S}
\]
in the derived category of $\OO_{X'_\bullet}$-modules such that 
$\mathcal{H}^i(\varphi)=C^{-1}$ for $i<p$.
\end{theorem}
\begin{proof}
To prove this theorem we simply check that all of the morphisms in the proof of Deligne-Illusie extend to 
morphisms on the simplicial level (see Section $\ref{subsec:review}$ for an outline of Deligne-Illusie 
and relevant notation).\\
\\
Let $\rho:X_n\rightarrow X_m$ be a face or a degeneracy map of $X_\bullet$.  To ease notation, we denote $X_n$ 
by $Y$ and $X_m$ by $X$.  In addition, we use $F$ to denote all relative Frobenii.\\
\\
To show that $\varphi^0$ extends to a morphism $\OO_{X'_\bullet}\rightarrow F_*\XBS^\bullet$, we show
\[
\xymatrix{
\rho'^*\mathcal{O}_{X'}\ar[r]^-{\rho'^*C^{-1}}\ar[d] & 
\rho'^*\mathcal{H}^0F_*\Omega^\bullet_{X/S} \ar[r]\ar[d] & 
\rho'^*F_*\Omega^\bullet_{X/S}\ar[d]\\
\mathcal{O}_{Y'}\ar[r]^-{C^{-1}} & \mathcal{H}^0F_*\Omega^\bullet_{Y/S} \ar[r] & 
F_*\Omega^\bullet_{Y/S}
}
\]
commutes.  The left square commutes by Lemma $\ref{l:Cartier}$.  The 
right square commutes since for any morphism 
$A^\bullet\rightarrow B^\bullet$ of complexes concentrated in non-negative degrees, the 
following diagram commutes
\[
\xymatrix{
\ker d^0_A \ar[r]\ar[d] & A^0\ar[d]\\
\ker d^0_B \ar[r] & B^0
}
\]
To show that $\varphi^i$ extends to a morphism on the simplicial level for $i>0$, we must check 
\[
\xymatrix{
\rho'^* \Omega^i_{X'/S}[-i]\ar[r]^-{\rho'^*a[-i]}\ar[d] & 
(\Omega^1_{X'/S})^{\otimes i}[-i]\ar[rr]^-{\rho'^*(\varphi^1)^{\otimes i}}\ar[d] & &
(F_*\Omega^\bullet_{X/S})^{\otimes i}\ar[r]^-{\rho'^*b}\ar[d] & 
F_*\Omega^\bullet_{X/S}\ar[d]\\
\Omega^i_{Y'/S}[-i]\ar[r]^-{a[-i]} & 
(\Omega^1_{Y'/S})^{\otimes i}[-i]\ar[rr]^-{(\varphi^1)^{\otimes i}} & &
(F_*\Omega^\bullet_{Y/S})^{\otimes i}\ar[r]^-{b} & F_*\Omega^\bullet_{Y/S}
}
\]
commutes.  It is clear that the outermost squares commute, and so we are reduced to 
checking the commutativity of 
\[
\xymatrix{
\rho'^*\Omega^1_{X'/S}[-1]\ar[r]^-{\rho'^*\varphi^1}\ar[d] & 
\rho'^*F_*\Omega^\bullet_{X/S}\ar[d]\\
\Omega^1_{Y'/S}[-1]\ar[r]^-{\varphi^1} & 
F_*\Omega^\bullet_{Y/S}
}
\]
Suppose now that Frobenius (for the simplicial scheme) lifts.  So, we have a commutative square 
\[
\xymatrix{
\tilde{Y}\ar[r]^-{\tilde{F}}\ar[d]_{\tilde{\rho}} & \tilde{Y'}\ar[d]^{\tilde{\rho}'}\\
\tilde{X}\ar[r]^-{\tilde{F}} & \tilde{X'}
}
\]
of $\tilde{S}$-schemes which pulls back to 
\[
\xymatrix{
Y\ar[r]^-{F}\ar[d]_{\rho} & Y'\ar[d]^{\rho'}\\
X\ar[r]^-{F} & X'
}
\]
over $S$.  In this case $\varphi^1=f$, and to check that it defines a morphism 
$\Omega^1_{X'_\bullet/S}[-1]\rightarrow F_*\Omega^1_{X_\bullet/S}$, we need to check that 
\[
\xymatrix{
\rho'^*\Omega^1_{X'/S}[-1]\ar[r]^-{\rho'^*f}\ar[d] & 
\rho'^*F_*\Omega^1_{X/S}\ar[d]\\
\Omega^1_{Y'/S}[-1]\ar[r]^-{f} & 
F_*\Omega^1_{Y/S}
}
\]
commutes.  Under these morphisms,
\[
\xymatrix{
dx_0\otimes1\ar@{|->}[r]\ar@{|->}[d] & x_0^{p-1}dx_0+du(x)\ar@{|->}[d]\\
d\rho(x_0)\otimes1 & \rho(x_0)^{p-1}d\rho(x_0)+d\rho(u(x))
}
\]
We see $d\rho(x_0)\otimes1$ is sent to  $\rho(x_0)^{p-1}d\rho(x_0)+d\rho(u(x))$ since 
\[
\tilde{F}^*(\til{\rho}(x)\otimes1)=\tilde{F}\tilde{\rho}'(x\otimes1)=
\tilde{\rho}\tilde{F}^*(x\otimes1)=\til{\rho}(x^p+\mathbf{p}(u(x)))=\til{\rho}(x)^p+\mathbf{p}(u(\til{\rho}(x))).
\]
Given two different choices $\til{F_1}$ and $\til{F_2}$ of $F$, we obtain a 
homotopy $h_{12}$ relating $f_1$ and $f_2$.  It is clear that $h_{12}$ extends to a morphism on the simplicial level 
since $h_{12}(dx_0\otimes1)=u_2(x)-u_1(x)$ and $\til{\rho}(\mathbf{p}(u_i(x)))=\mathbf{p}(u_i(\til{\rho}(x)))$.\\
\\
We now need to handle the general case when Frobenius does not lift.  We begin by proving that Frobenius lifts \'etale 
locally on $X_\bullet$.  To do so, we can, by Proposition \ref{l:Artinfroblift}, assume that there is a lift $\til{F}$ of 
$F_{\mf{X}/S}$.  Let $\mathcal{U}_0=\{U_i\}$ be a Zariski cover of $X_0$ where Frobenius lifts and let $\til{F}_i$ be 
a lift of $F_{U_i/S}$.  Then $\mathcal{U}_n:=\{U_{i_1}\times_{\mf{X}}\dots\times_{\mf{X}} U_{i_n}\}$ is a 
Zariski cover of $X_n$ and $\til{F}_{i_1}\times_{\til{F}}\dots\times_{\til{F}}\til{F}_{i_n}$ is a lift of Frobenius 
on $U_{i_1}\times_{\mf{X}}\dots\times_{\mf{X}} U_{i_n}$.  Moreover, these lifts of Frobenius are 
compatible so we see that Frobenius for the simplicial scheme does lift \'etale locally.\\
\\
To finish the proof of the theorem, we need only prove the commutativity of 
\[
\xymatrix{
\rho'^*F_*\Omega_{X/S}^\bullet \ar[r]^-{\simeq}\ar[d] & 
\textrm{Tot}(\rho'^*F_*\check{\mathcal{C}}^\bullet(\mathcal{U}_m,\Omega^\bullet_{X/S}))\ar[d] & 
\rho'^*\Omega^1_{X'/S}[-1]\ar[l]_-{\rho'^*\Phi}\ar[d]\\
F_*\Omega_{Y/S}^\bullet \ar[r]^-{\simeq} & 
\textrm{Tot}(F_*\check{\mathcal{C}}^\bullet(\mathcal{U}_n,\Omega^\bullet_{Y/S})) & 
\Omega^1_{Y'/S}[-1]\ar[l]_-{\Phi}
}
\]
The right square commutes because the 
$\Phi$ are defined in terms of the $f$'s and $h$'s.  The middle vertical map 
is induced by the morphism of the respective double complexes given by 
\[
(\omega_{1,s}\wedge\dots\wedge\omega_{a,s})_{s\in S_m}\mapsto
(\rho_s^{st}(\omega_{1,s}\wedge\dots\wedge\rho_s^{st}(\omega_{a,s}))_{s\in S_m,t\in S_n}
\]
where $\rho_s^{st}:U_s\times_{\mf{X}} U_t\rightarrow U_s$ and $S_k$ is the symmetric group.  So, under 
the morphisms in the left square, 
\[
\xymatrix{
\omega_1\wedge\dots\omega_a \ar@{|->}[r]\ar@{|->}[d] & 
(\omega_1|U_s\wedge\dots\wedge\omega_a|U_s)_s\\
\rho(\omega_1)\wedge\dots\rho(\omega_a) \ar@{|->}[r] & 
(\rho(\omega_1)|U_s\times U_t\wedge\dots\wedge\rho(\omega_a)|U_s\times U_t)_{s,t}
}
\]
Under the middle vertical map, $(\omega_1|U_s\wedge\dots\wedge\omega_a|U_s)_s$ is 
sent to 
$(\rho_s^{st}(\omega_1|U_s)\wedge\dots\wedge\rho_s^{st}(\omega_a|U_s))_{s,t}$.  But 
\[
\xymatrix{
U_s\times_{\mf{X}} U_t \ar[r]\ar[d]_{\rho_s^{st}} & Y\ar[d]^{\rho}\\
U_s\ar[r] & X
}
\]
commutes, so this completes the proof.
\end{proof}
%\begin{remark}\label{rmk:liftablecover}\emph{If $\mf{X}$ is a smooth Artin stack which lifts mod $p^2$, then there automatically exists a smooth cover $X\rightarrow\mf{X}$ by a smooth scheme such that the cover lifts mod $p^2$.  This can be seen as follows.  Let $Y\rightarrow\mf{X}$ be any smooth cover by a smooth scheme $Y$ and let $\bigcup U_i=Y$ be a Zariski cover of $Y$ by open affine subschemes.  We can then take $X=\coprod U_i$.}\end{remark}

\section{de Rham Theory for Tame Stacks}
%Given an Artin stack $\mf{X}$ over a scheme $S$, the lisse-\'etale site $Lis-et(\mf{X})$ of $\mf{X}$ is defined as follows.  
%The underlying category is the full subcategory of $\mf{X}$-schemes which are smooth over $\mf{X}$, and coverings are given 
%by $\{U_i\rightarrow U\}$ where each $U_i\rightarrow U$ is surjective and $\coprod U_i\rightarrow U$ is surjective.  
%As shown in \cite[12.2.1]{lmb}, giving a sheaf on the lisse-\'etale site of $\mf{X}$ is equivalent to giving the data of a 
%sheaf $\mathcal{F}_U\in U_{et}$ for all $U\in Lis-et(\mf{X})$, and transition maps 
%\[
%\phi_f:f^{-1}\mathcal{F}_U\rightarrow\mathcal{F}_V
%\]
%for all morphisms $f:V\rightarrow U$ in $Lis-et(\mf{X})$ such that $\phi_f$ is an isomorphism if $f$ is \'etale and 
%$\phi_g\circ g^*\phi_f=\phi_{gf}$ for all composable morphisms $f$ and $g$.
Let $S$ be a scheme and $\mf{X}\rightarrow\mathcal{Y}$ a morphism of Artin stacks over $S$.  We denote by 
$\varpi_{\mf{X}/\mathcal{Y}}^1$ the sheaf of $\OO_\mf{X}$-modules on the lisse-\'etale site of $\mf{X}$ such that 
$\varpi_{\mf{X}/\mathcal{Y}}^1|U_{et}=\Omega^1_{U/\mathcal{Y}}$ for all $U$ smooth over $\mf{X}$.  We define 
$\varpi_{\mf{X}/\mathcal{Y}}^s$ to be $\bigwedge^s\varpi_{\mf{X}/\mathcal{Y}}^1$.  Given a morphism 
$f:V\rightarrow U$ of smooth $\mf{X}$-schemes, note that the transition function 
\[
f^*\Omega^1_{U/\mathcal{Y}}\longrightarrow\Omega^1_{V/\mathcal{Y}}
\]
need not be an isomorphism, and so the $\varpi_{\mf{X}/\mathcal{Y}}^s$ are never coherent.  Note also that $\XX^1$ 
is not the zero sheaf.\\
\\
As mentioned in the introduction, the sheaf $\XS^1$ gives us a naive de Rham complex $\XS^\bullet$.  In this section 
we prove that when $S$ is spectrum of a perfect field of characteristic $p$, the hypercohomology spectral sequence 
\[
E_1^{st}=H^t(\XS^s)\Rightarrow H^n(\XS^\bullet)
\]
degenerates for smooth proper tame stacks $\mf{X}$ 
that lift mod $p^2$.  The reason the proof of Corollary \ref{cor:dmdegen} and the Deligne-Illusie result proved in the 
last section do not immediately imply the degeneracy of this spectral sequence 
%for all smooth proper \emph{Artin} stacks satisfying the hypotheses of Theorem \ref{thm:disimp}, let alone tame stacks, 
is that, as mentioned above, the $\XS^s$ are not coherent, and so we do not yet know that the $H^t(\XS^s)$ and 
$H^n(\XS^\bullet)$ are finite-dimensional $k$-vector spaces.  The main goal of this section, which implies the degeneracy of 
the above spectral sequence, is to prove that they are 
%When $\mf{X}$ is smooth and tame, however, we show that 
by comparing them with the cohomology of the cotangent complex.\\
\\
We begin by proving three general lemmas and a corollary which require no assumptions on the base scheme $S$.  The first two 
lemmas are concerned with relative cohomological descent.  For background material on cohomological descent, we refer the 
reader to \cite[\S 2]{sheaves} and \cite[\S 6]{conrad}.\\
\\
In what follows, given a smooth hypercover $a:X_\bullet \rightarrow\mf{X}$ of an Artin stack by a simplicial algebraic space, 
$\mf{X}_{lis-et}|X_s$ denotes the topos of sheaves over the representable sheaf defined by $X_s$ and 
$\mf{X}_{lis-et}|X_\bullet$ denotes the associated simplicial topos.
\begin{lemma}
\label{l:first}
Let $\mf{X}$ be an Artin stack over $S$ and let $a:X_\bullet \rightarrow\mf{X}$ be a smooth hypercover by a simplicial 
algebraic space.  
If $f:\mf{X}\rightarrow M$ is a morphism to a scheme, then for any $\mathcal{F}_\bullet\in Ab(\mf{X}_{lis-et}|X_\bullet)$, 
there is a spectral sequence
\[
E_1^{st}=R^t(fa_s)_*(\mathcal{F}_s|X_{s,et})\Rightarrow \epsilon_*R^n(f_* a_*)\mathcal{F}_\bullet,
\]
where $\epsilon:M_{lis-et}\rightarrow M_{et}$ is the canonical morphism of topoi.  If 
$\mathcal{F}_\bullet=a^*\mathcal{F}$ for some $\mathcal{F}\in Ab(\mf{X}_{lis-et})$, then 
$\epsilon_*R^n(f_* a_*)\mathcal{F}_\bullet=\epsilon_*R^nf_*\mathcal{F}$.
\end{lemma}
\begin{proof}
Let $\eta_s:\mf{X}_{lis-et}|X_s\rightarrow X_{s,et}$ be the canonical morphism of topoi and note that 
\[
\xymatrix{
Ab(X_{s,et})\ar^{(fa_s)_*}[drrr] & & & \\
Ab(\mf{X}_{lis-et}|X_s)\ar^{{\eta_s}_*}[u]\ar_-{{a_s}_*}[r] & Ab(\mf{X}_{lis-et})\ar_{f_*}[r] & 
Ab(M_{lis-et})\ar_-{\epsilon_*}[r] & Ab(M_{et})
}
\]
commutes.  By general principles (see proof of \cite[Thm 6.11]{conrad}), there is a spectral sequence
\[
E_1^{st}=R^t(\epsilon_*f_*a_{s*})(\mathcal{F}_s)\Rightarrow R^n(\epsilon_*f_* a_*)\mathcal{F}_\bullet.
\]
As $\epsilon_*$ is exact, $R^n(\epsilon_*f_* a_*)\mathcal{F}_\bullet=\epsilon_*R^n(f_* a_*)\mathcal{F}_\bullet$.  
Since ${\eta_s}_*$ is exact and takes injectives to injectives, the commutativity of the above diagram implies that 
$E_1^{st}\simeq R^t(fa_s)_*(\mathcal{F}_s|X_{s,et})$, which shows the existence of our desired spectral sequence.  Lastly, 
since
\[
a^*:Ab(\mf{X}_{lis-et})\rightarrow Ab(\mf{X}_{lis-et}|X_\bullet)
\]
is fully faithful, it follows (\cite[Lemma 6.8]{conrad}) that $Ra_*a^*=id$.  As a result, 
$\epsilon_* R(f_*a_*)a^*\mathcal{F}=\epsilon_* Rf_*\mathcal{F}$.
\end{proof}
\begin{lemma}
\label{problem}
With notation and hypotheses as in Lemma \ref{l:first}, we have
\[
R^n(fa)_*(\eta_*\mathcal{F}_\bullet)=\epsilon_*R^n(f_*a_*)\mathcal{F}_\bullet,
\]
where $\eta:\mf{X}_{lis-et}|X_\bullet \rightarrow X_{\bullet,et}$ is the canonical morphism of topoi.
\end{lemma}
\begin{proof}
We see that the diagram 
\[
\xymatrix{
Ab(X_{\bullet,et})\ar^{(fa)_*}[drrr] & & & \\
Ab(\mf{X}_{lis-et}|X_\bullet)\ar^{\eta_*}[u]\ar_-{a_*}[r] & Ab(\mf{X}_{lis-et})\ar_{f_*}[r] & 
Ab(M_{lis-et})\ar_-{\epsilon_*}[r] & Ab(M_{et})
}
\]
commutes.  It follows that 
\[
R(fa)_*(\eta_*\mathcal{F}_\bullet)=R(fa)_*(R\eta_*\mathcal{F}_\bullet)=\epsilon_*R(f_*a_*)\mathcal{F}_\bullet,
\]
as $\epsilon_*$ and $\eta_*$ are exact and take injectives to injectives.
%\\
%\\
%Alternatively, $(f\circ a)_*:Ab(X_{\bullet,et})\rightarrow Ab(M_{et})$ is an augmented simplicial topos, 
%so we have a spectral sequence $E_1^{pq}=R^q(f_* (a_p)_*)(\mathcal{F}_p)
%\Rightarrow R^n(f_* a_*)(\mathcal{F}_\bullet)$.
\end{proof}
Using Lemma \ref{l:first}, we prove a base change result for sheaves on an Artin stack which are not necessarily 
quasi-coherent, but are level-by-level quasi-coherent on a smooth hypercover of the stack.
\begin{corollary}
\label{cor:basechange}
Let $f:\mf{X}\rightarrow M$ be a morphism from an Artin stack to a scheme and let $a:X_\bullet\rightarrow \mf{X}$ be a smooth 
hypercover by a simplicial algebraic space.  Let $h:T\rightarrow M$ be an \'etale morphism and consider the diagram
\[
\xymatrix{
Y_\bullet\ar^{j}[r]\ar_{b}[d] & X_\bullet\ar^{a}[d]\\
\mathcal{Y}\ar[r]^{i}\ar_{g}[d] & \mf{X}\ar^{f}[d]\\
T\ar[r]^h & M
}
\]
where all squares are cartesian.  If $\mathcal{F}$ is an $\OO_{\mf{X}}$-module such that each $\mathcal{F}|X_{s,et}$ is 
quasi-coherent, then the canonical map 
\[
h^*\epsilon_*R^nf_*\mathcal{F}\longrightarrow \alpha_*R^ng_*i^*\mathcal{F}
\]
is an isomorphism, where $\epsilon$ and $\alpha$ denote the canonical morphisms of topoi $M_{lis-et}\rightarrow M_{et}$ and 
$T_{lis-et}\rightarrow T_{et}$, respectively.
\end{corollary}
\begin{proof}
By Lemma \ref{l:first}, we have a spectral sequence 
\[
E_1^{st}=R^t(fa_s)_*(\mathcal{F}|X_{s,et})\Rightarrow \epsilon_*R^nf_*\mathcal{F}.
\]
Applying $h^*$, we obtain another spectral sequence
\[
{}^{\backprime} E_1^{st}=h^*R^t(fa_s)_*(\mathcal{F}|X_{s,et})\Rightarrow h^*\epsilon_*R^nf_*\mathcal{F}.
\]
Flat base change shows 
\[
h^*R^t(fa_s)_*(\mathcal{F}|X_{s,et})=R^t(gb_s)_*j_s^*(\mathcal{F}|X_{s,et})=R^t(gb_s)_*(i^*\mathcal{F})|Y_{s,et}.
\]
Another application of Lemma \ref{l:first} then shows that ${}^{\backprime} E$ in fact abuts to $\alpha_*R^ng_*i^*\mathcal{F}$.
\end{proof}
Before stating the next lemma, we introduce the following definitions.  Let $Z$ be an $S$-scheme equipped with an action 
$\rho:G\times_S Z\rightarrow Z$ of a smooth linearly reductive group scheme $G$ over $S$ and let $p:G\times Z\rightarrow Z$ be the 
projection.  We 
denote by ($G$-$lin\ \OO_{Z_{et}}$-$mod$) the category of $G$-linearized $\OO_{Z_{et}}$-modules.  That is, the category 
of quasi-coherent $\OO_{Z_{et}}$-modules $\mathcal{F}$ together with an isomorphism 
$\phi:p^*\mathcal{F}\rightarrow\rho^*\mathcal{F}$ satisfying a cocycle condition.  From such a $\phi$ we can define a 
``coaction map''
\[
\sigma:\mathcal{F}\longrightarrow p_*p^*\mathcal{F} \stackrel{p_*(\phi)}{\longrightarrow} p_*\rho^*\mathcal{F}=
\mathcal{F}\otimes_{\OO_Z,\rho}\OO_{Z\times_S G}
\]
which satisfies an associativity relation as in \cite[p.31]{git}; here $\mathcal{F}\rightarrow p_*p^*\mathcal{F}$ 
is the canonical map.  Letting $f:Z\rightarrow Z/G$ be the natural map, we define the $G$-invariants $\mathcal{F}^G$ 
of $\mathcal{F}$ to be the equalizer of
\[
f_*\sigma:f_*\mathcal{F}\longrightarrow f_*p_*\rho^*\mathcal{F}=f_*p_*p^*\mathcal{F}
\]
and $f_*$ of the canonical map $s\mapsto s\otimes1$.\\
\\
If $Y$ is also an $S$-scheme equipped with a $G$-action and $h:Z\rightarrow Y$ is a $G$-equivariant map over $S$, then for 
every $G$-linearized $\OO_{Z_{et}}$-module $\mathcal{F}$, there is a natural $G$-linearization on $h_*\mathcal{F}$.  So, 
we have a commutative diagram of categories 
\[
\xymatrix{
(G\textrm{-}lin\ \OO_{Z_{et}}\textrm{-}mod) \ar[r]\ar[d]_{h_*} & (\OO_{Z_{et}}\textrm{-}mod) \ar[d]^{h_*}\\
(G\textrm{-}lin\ \OO_{Y_{et}}\textrm{-}mod) \ar[r] & (\OO_{Y_{et}}\textrm{-}mod)
}
\]
where the horizontal arrows are the obvious forgetful functors.  If $g:Z/G\rightarrow Y/G$ denotes the map induced by $h$, 
then it is not hard to see that $(h_*\mathcal{F})^G=g_*\mathcal{F}^G$.  In particular, $\mathcal{F}^G=(f_*\mathcal{F})^G$ 
where $G$ acts trivially on $Z/G$.  Note that for any sheaf $\mathcal{G}$ of $\OO_{Z/G}$-modules, 
$f^*\mathcal{G}$ comes equipped with a canonical $G$-linearization.  If the $G$-action on $Z$ is free, so that $f$ is a 
$G$-torsor, then $(f^*\mathcal{G})^G=(f_*f^*\mathcal{G})^G=\mathcal{G}$.\\
\\
By descent theory, ($G$-$lin\ \OO_{Z_{et}}$-$mod$) is equivalent to the category of quasi-coherent sheaves on $[Z/G]$.  Under 
this equivalence, taking $G$-invariants in the above sense corresponds to pushing forward to the coarse space $Z/G$.\\
%As a result, if $G$ is linearly reductive, then $R^n(-)^G(\mathcal{F})=0$ for $n>0$.\\
\\
When the action of $G$ 
on $Z$ is trivial, we denote ($G$-$lin\ \OO_{Z_{et}}$-$mod$) by ($G$-$\OO_{Z_{et}}$-$mod$).  We can similarly define the 
categories ($G$-$lin\ \OO_{Z_{\bullet,et}}$-$mod$) and ($G$-$\OO_{Z_{\bullet,et}}$-$mod$) for simplicial schemes $Z_\bullet$.
\begin{lemma}
\label{Martin}
Let $U$ be a smooth $S$-scheme with an action of a smooth affine linearly reductive group scheme $G$ over $S$.  Let 
$\mf{X}=[U/G]$ and $a:X_\bullet\rightarrow \mf{X}$ be the hypercover obtained by taking fiber products of $U$ over $\mf{X}$.  
Consider the diagram
\[
\xymatrix{
Y_\bullet\ar^{\pi}[r]\ar_{b}[d] & X_\bullet\ar^{a}[d]\\
U\ar[r]^{a_0}\ar_{g}[dr] & \mf{X}\ar^{f}[d]\\
 & M
}
\]
where the square is cartesian and $M$ is a scheme.  Then
\[
R^n(fa)_*\mathcal{F}_\bullet=(R^n(gb)_*\pi^*\mathcal{F}_\bullet)^G
\]
for all $\OO_{X_{\bullet,et}}$-modules $\mathcal{F}_\bullet$ such that the $\mathcal{F}_s$ are quasi-coherent.  
\end{lemma}
\begin{proof}
Note that the following diagram 
\[
\xymatrix{
(G\textrm{-}lin\ \OO_{Y_{\bullet,et}}\textrm{-}mod) \ar[rr]^-{(gb)_*} \ar[d]_{\pi_*} & & 
(G\textrm{-}\OO_{M_{et}}\textrm{-}mod) \ar[d]^{(-)^G}\\
(G\textrm{-}\OO_{X_{\bullet,et}}\textrm{-}mod)\ar[r]^-{(-)^G} & (\OO_{X_{\bullet,et}}\textrm{-}mod) \ar[r]^-{(fa)_*}  & 
(\OO_{M_{et}}\textrm{-}mod)
}
\]
of categories commutes.  As a result, 
\[
R(fa)_*R(-)^G(R\pi_*\pi^*\mathcal{F}_\bullet)=R(-)^G(R(gb)_*\pi^*\mathcal{F}_\bullet)=(R(gb)_*\pi^*\mathcal{F}_\bullet)^G,
\]
where the second equality holds because $R(gb)_*\pi^*\mathcal{F}_\bullet$ has quasi-coherent cohomology.  It suffices then 
to prove 
\[
\mathcal{F}_\bullet=R(-)^G(R\pi_*\pi^*\mathcal{F}_\bullet).
\]
We begin by showing $R\pi_*\pi^*\mathcal{F}_\bullet=\pi_*\pi^*\mathcal{F}_\bullet$.  Let
\[
0\rightarrow\pi^*\mathcal{F}_\bullet\rightarrow \mathcal{I}_\bullet^0\rightarrow \mathcal{I}_\bullet^1\rightarrow\dots
\]
be an injective resolution of $\OO_{X_{\bullet,et}}$-modules.  To show 
$R^n\pi_*\pi^*\mathcal{F}_\bullet=0$ for $n>0$, we need only do so after restricting to each 
level $X_s$.  Since the restriction functor $res_s:Ab(X_{\bullet,et})\rightarrow Ab(X_{s,et})$ is exact, we see
\[
res_s R^n\pi_*\pi^*\mathcal{F}_\bullet= res_s \HH^n(\pi_*\mathcal{I}_\bullet^\bullet)= 
\HH^n(\pi_*\mathcal{I}_s^\bullet)= R^n\pi_*\pi^*\mathcal{F}_s=0,
\]
where the last equality holds because $\pi$ is affine and $\mathcal{F}_s$ is quasi-coherent.\\
\\
A similar argument shows $R(-)^G(\pi_*\pi^*\mathcal{F}_\bullet)=(\pi_*\pi^*\mathcal{F}_\bullet)^G$ as every 
$\pi_*\pi^*\mathcal{F}_s$ is quasi-coherent.  The lemma then follows from the fact that $\pi$ is a $G$-torsor, and so 
$(\pi_*\pi^*\mathcal{F}_\bullet)^G=\mathcal{F}_\bullet$.
\end{proof}
For the rest of the section, we let $S=\Spec k$, where $k$ is a perfect field of characteristic $p$.  
We remind the reader that if $\mf{X}$ is a smooth Artin stack and $X_\bullet\rightarrow\mf{X}$ is a hypercover, then 
the cotangent complex $L_{\mf{X}/S}$ of the stack (\cite[\S 8]{sheaves}) is the bounded complex of $\OO_{\mf{X}}$-modules 
with quasi-coherent cohomology such that 
\[
L_{\mf{X}/S}|X_{\bullet,et}=\Omega^1_{X_\bullet/S}\rightarrow\Omega^1_{X_\bullet/\mf{X}}
\]
with $\Omega^1_{X_\bullet/S}$ in degree $0$; that is,
\[
L_{\mf{X}/S}=\XS^1\rightarrow \XX^1.
\]
In Theorem \ref{main2}
below, we compare $\XS^s$ with $\pLXS$, the $s^{th}$ derived exterior power of $\LXS$.  Given an 
abelian category $\mathcal{A}$, the derived exterior powers $L\bigwedge^s$, as well as the derived symmetric powers $LS^s$, 
of a complex $E\in D^-(\mathcal{A})$ are defined in \cite[I.4.2.2.6]{cotangent}.  Since $\LXS$ is not concentrated in 
negative degrees, we cannot directly define $\pLXS$; however, it is shown in \cite[I.4.3.2.1]{cotangent} that for 
$E\in D^-(\mathcal{A})$, 
\[
LS^s(E[1])=(L\bigwedge^s E)[s]
\]
so we may define $\pLXS$ as $LS^s(\LXS[1])[-s]$.  It follows, then, from \cite[I.4.3.1.7]{cotangent} that 
\[
\pLXS=\XS^s\rightarrow \XS^{s-1}\otimes\XX^1\rightarrow\dots\rightarrow \XS^1\otimes S^{s-1}\XX^1\rightarrow S^s\XX^1
\]
with $\XS^s$ in degree $0$.  Note that we have a canonical map from $\pLXS$ to $\XS^s$.\\
\\
We remark that $\pLXS\in D^b_{coh}(\mf{X})$ for all $s$.  This can be seen as follows.  We have an exact triangle
\[
a_0^*\LXS\longrightarrow L_{X_0/S}\longrightarrow \Omega^1_{X/\mf{X}}.
\]
By \cite[II.2.3.7]{cotangent}, $L_{X_0/S}$ has coherent cohomology.  Since $\Omega^1_{X/\mf{X}}$ is coherent and coherence can 
be checked smooth locally, we see $L_{X/S}$ and hence all $\pLXS$ are in $D^b_{coh}(\mf{X})$.\\
\\
We are now ready to prove the comparison theorem.
\begin{theorem}
\label{main2}
If $\mf{X}$ is smooth and tame over $S$ and $f:\mf{X}\rightarrow M$ is its coarse space, then the canonical map
\[
\epsilon_*R^tf_*(\pLXS)\longrightarrow\epsilon_*R^tf_*\XS^s
\]
is an isomorphism.
\end{theorem}
\begin{proof}
By \cite[Prop 5.2]{cstlr}, 
%Lemma \ref{l:semidirect}, 
there exists an \'etale cover $h:T\rightarrow M$ and a cartesian diagram
\[
\xymatrix{
[U/G]\ar[r]\ar_{g}[d] & \mf{X}\ar^{f}[d]\\
T\ar[r]^{h} & M
}
\]
where $G$ is linearly reductive, affine, and smooth over $S$.  Since $\mf{X}$ and $G$ are smooth, we see that $U$ is as 
well.  Let $\mathcal{Y}=[U/G]$ and let $\varphi:\epsilon_*R^tf_*(\pLXS)\rightarrow \epsilon_*R^tf_*\XS^s$ be the 
canonical map.  By Corollary \ref{cor:basechange}, we see that $h^*\varphi$ is the canonical map 
\[
\epsilon_*R^tg_*(\bigwedge^sL_{\mathcal{Y}/S})\longrightarrow\epsilon_*R^tg_*\varpi_{\mathcal{Y}/S}^s.
\]
To show that $\varphi$ is an isomorphism, we can therefore assume $\mf{X}=[U/G]$ and $M=T$.\\
\\
To prove the theorem, it suffices 
to show $\epsilon_*R^tf_*(\XS^{s-k}\otimes S^k\XX^1)=0$ for all $k>0$ and all $t$.  With notation as in 
Lemma \ref{Martin}, we see 
\[
\epsilon_*R^tf_*(\XS^{s-k}\otimes S^k\XX^1)=\epsilon_*R^t(f_*a_*)a^*(\XS^{s-k}\otimes S^k\XX^1)=
R^t(fa)_*(\Omega^{s-k}_{X_\bullet/S}\otimes S^k\Omega^1_{X_\bullet/\mf{X}}),
\]
where the first equality is by Lemma \ref{l:first} and the second is by Lemma \ref{problem}.  It now follows from Lemma 
\ref{Martin} that 
\[
R^t(fa)_*(\Omega^{s-k}_{X_\bullet/S}\otimes S^k\Omega^1_{X_\bullet/\mf{X}})=
(R^t(gb)_*(\pi^*\Omega^{s-k}_{X_\bullet/S}\otimes S^k\Omega^1_{Y_\bullet/U}))^G.
\]
Fix $t$ and $k>0$.  It suffices then to prove by (strong) induction on $s$ that for every flat $\OO_{X_\bullet}$-module 
$\mathcal{G}$ which is restriction to $Y_{\bullet,et}$ of some $\OO$-module $\mathcal{F}$ on the lisse-\'etale site of $U$,
\[
R^n(gb)_*(\pi^*\Omega^s_{X_\bullet/S}\otimes\mathcal{G}\otimes S^k\Omega^1_{Y_\bullet/U})=0.
\]
We begin with the case $s=0$, which is handled separately.  An application of Lemmas \ref{l:first} and \ref{problem} shows
\[
R^n(gb)_*(\mathcal{G}\otimes S^k\Omega^1_{Y_\bullet/U})=
\epsilon_*R^ng_*(\mathcal{F}\otimes S^k\varpi^1_{U/U}).
\]
If we let $\alpha:U_{lis-et}\rightarrow U_{et}$ be the canonical morphism of topoi, we see then that
\[
\epsilon_*R^ng_*(\mathcal{F}\otimes S^k\varpi^1_{U/U})=R^ng_*(\alpha_*\mathcal{F}\otimes S^k\Omega^1_{U/U})=0,
\]
where the last equality holds since $k>0$.\\
\\
Assume now that $s>0$.  Since $\pi$ is smooth, we have a short exact sequence 
\[
0 \longrightarrow\pi^*\Omega^1_{X_\bullet/S} \longrightarrow\Omega^1_{Y_\bullet/S} 
\longrightarrow\Omega^1_{Y_\bullet/X_\bullet} \longrightarrow 0.
\]
As a result, we have a filtration $\Omega^s_{Y_\bullet/S}\supset \mathcal{K}^1\supset\dots\supset\mathcal{K}^s\supset 0$ 
with $\mathcal{K}^s=\pi^*\Omega^s_{X_\bullet/S}$ and short exact sequences 
\[
0 \longrightarrow \mathcal{K}^1 \longrightarrow\Omega^s_{Y_\bullet/S} 
\longrightarrow\Omega^s_{Y_\bullet/X_\bullet} \longrightarrow 0
\]
\[
0 \longrightarrow \mathcal{K}^2 \longrightarrow \mathcal{K}^1
\longrightarrow \pi^*\Omega^1_{X_\bullet/S}\otimes\Omega^{s-1}_{Y_\bullet/X_\bullet} \longrightarrow 0
\]
\[
\vdots
\]
\[
0 \longrightarrow \pi^*\Omega^s_{X_\bullet/S} \longrightarrow \mathcal{K}^{s-1}
\longrightarrow \pi^*\Omega^{s-1}_{X_\bullet/S}\otimes\Omega^1_{Y_\bullet/X_\bullet} \longrightarrow 0.
\]
Since $\mathcal{G}\otimes S^k\Omega^1_{Y_\bullet/U}$ is flat, tensoring each of the above short exact sequences by it results 
in a new list of short exact sequences.  Since 
\[
\Omega^1_{Y_\bullet/S}\otimes\mathcal{G}\otimes S^k\Omega^1_{Y_\bullet/U}=
(\varpi^1_{U/S}\otimes\mathcal{F}\otimes S^k\varpi^1_{U/U})|Y_{\bullet,et}
\]
and
\[
\Omega^1_{Y_\bullet/X_\bullet}\otimes\mathcal{G}\otimes S^k\Omega^1_{Y_\bullet/U}=
(L_{U/\mf{X}}\otimes\mathcal{F}\otimes S^k\varpi^1_{U/U})|Y_{\bullet,et},
\]
the $s=0$ case shows 
\[
R^n(gb)_*(\Omega^1_{Y_\bullet/S}\otimes\mathcal{G}\otimes S^k\Omega^1_{Y_\bullet/U})=
R^n(gb)_*(\Omega^1_{Y_\bullet/X_\bullet}\otimes\mathcal{G}\otimes S^k\Omega^1_{Y_\bullet/U})=0.
\]
As a result, $R^n(gb)_*(\mathcal{K}^1\otimes\mathcal{G}\otimes S^k\Omega^1_{Y_\bullet/U})=0$.  Using the inductive 
hypothesis, we conclude
\[
R^n(gb)_*(\mathcal{K}^i\otimes\mathcal{G}\otimes S^k\Omega^1_{Y_\bullet/U})=0
\]
for all $i$, in particular for $i=s$.
\end{proof}
\begin{corollary}
\label{main}
If $\mf{X}$ is a smooth proper tame stack over $S$, then $H^t(\XS^s)$ and $H^n(\XS^\bullet)$ are 
finite-dimensional $k$-vector spaces for all $s$, $t$, and $n$.
\end{corollary}
\begin{proof}
Let $f:\mf{X}\rightarrow M$ be the coarse space of $\mf{X}$.  For each $s$, there is a Leray spectral sequence 
\[
E_2^{ij}=H^i(\epsilon_*R^jf_*\XS^s)\Rightarrow H^t(\XS^s).
\]
By Theorem \ref{main2}, the canonical map
\[
\epsilon_*R^jf_*(\pLXS)\longrightarrow \epsilon_*R^jf_*\XS^s
\]
is an isomorphism.  As we remarked above, $\pLXS\in D^b_{coh}(\mf{X})$.  
Since $f$ is proper by Keel-Mori \cite{km}, and $M$ is proper by \cite[Prop 2.10]{hom}, we see the $E_2^{ij}$ are 
finite-dimensional $k$-vector spaces.  It follows that $H^t(\XS^s)$ 
is a finite-dimensional $k$-vector space for every $s$ and $t$.\\
\\
Since the morphisms in the complex $\XSB$ are $k$-linear, the hypercohomology spectral sequence 
\[
E_1^{st}=H^t(\XS^s)\Rightarrow H^n(\XS^\bullet)
\]
consists of finite-dimensional $k$-vector spaces with $k$-linear maps.  As a result, $H^n(\XS^\bullet)$ is a 
finite-dimensional $k$-vector space as well.
\end{proof}
\begin{theorem}
\label{cor:tamedegen}
Let $\mf{X}$ be a smooth proper tame stack over $S$ that lifts mod $p^2$.  Then the 
Hodge-de Rham spectral sequence
\[
E_1^{st}=H^t(\XS^s)\Rightarrow H^n(\XS^\bullet)
\]
degenerates for $s+t<p$.
\end{theorem}
\begin{proof}
  %By Remark~\ref{rmk:liftablecover}, there exists a smooth cover
  Choose a smooth cover 
  $X\rightarrow\mf{X}$ by a smooth scheme such that the cover lifts
  mod $p^2$.  Theorem~\ref{thm:disimp} now shows
\[
\bigoplus_{s<p}\Omega^s_{X'_\bullet/S}[-s]\simeq
\tau_{<p}F_*\Omega^\bullet_{X_\bullet/S},
\]
where $X_\bullet$ is obtained from $X$ by taking fiber products over
$\mf{X}$.  Since $H^t(\XS^s)=H^t(\XBS^s)$ and
$H^n(\XS^\bullet)=H^n(\XBS^\bullet)$, we see that for $n<p$,
\begin{multline*}
\dim_k H^n(\XS^\bullet)=\sum_{s+t=n} \dim_k H^t(\Omega^s_{X'_\bullet/S})=
\sum_{s+t=n} \dim_k H^t(\XBS^s)\\
=\sum_{s+t=n} \dim_k H^t(\XS^s),
\end{multline*}
which proves the degeneracy of the spectral sequence by 
Lemma~\ref{l:homological}.
\end{proof}

\section{de Rham Theory for Schemes\\ with Isolated Linearly Reductive Singularities}
Let $k$ be a perfect field of characteristic $p$ and let $S=\Spec k$.
%If $G$ is a finite linearly reductive group scheme over $k$, 
%then by \cite[Thm 2.19]{tame}, $G=\Delta\rtimes H$, where $\Delta$ is a finite diagonalizable group scheme and $H$ is a 
%finite \'etale constant group scheme.  Let $\alpha:\Delta\hookrightarrow \mathbb{G}_m^r$ be an embedding.  We define 
%$G^\alpha$ as the pushout of $\mathbb{G}_m^r\coprod_{alpha,\Delta}G$ so that $G^\alpha=\mathbb{G}_m^r \rtimes H$.
%\[
%\xymatrix{
%1\ar[r]\ar[d] & \Delta\ar[r]\ar^{\iota}[d] & G\ar[r]\ar[d] & H\ar[r]\ar^{id}[d] & 1\\
%1\ar[r] & \mathbb{G}_m^r \ar[r] & G^{\alpha}\ar[r] & H\ar[r] & 1
%}
%\]
\begin{definition}
\emph{We say a scheme $M$ over $S$ has $\emph{linearly}$ $\emph{reductive}$ $\emph{singularities}$ if 
there is an \'etale cover $\{U_i/G_i\rightarrow M\}$, where the $U_i$ are smooth over $S$ and the 
$G_i$ are linearly reductive group schemes which are finite over $S$.}
%and the $U_i\otimes^{G_i}G_i^{\alpha_i}$ are smooth 
%for all embeddings $\alpha_i$.}
\end{definition}
Note that if $M$ has linearly reductive singularities, then it is automatically normal and in fact Cohen-Macaulay by 
\cite[p.115]{cm}.\\
\\
Our goal in this section is to prove that if $M$ is proper over $S$, and $j:M^0\rightarrow M$ is its smooth locus, 
then under suitable liftability conditions, the hypercohomology spectral sequence 
$E_1^{st}=H^t(j_*\DR^s)\Rightarrow H^n(j_*\DRB)$ 
degenerates.  
\subsection{Relationship with Tame Stacks, and the Cartier Isomorphism}
We begin by recalling the relationship between tame stacks and schemes with linearly reductive singularities:
\begin{theorem}[{\cite[Thm 1.10]{cstlr}}]
\label{prop:cslrs}
Let $M$ be an $S$-scheme with linearly reductive singularities.  Then it is the coarse space of a smooth tame stack 
$\mf{X}$ such that $f^0$ in the diagram 
\[
\xymatrix{
\mf{X}^0\ar[r]^{j^0}\ar[d]_{f^0} & \mf{X}\ar[d]^f \\
M^0\ar[r]_j & M
}
\]
is an isomorphism, where $\mf{X}^0=M^0\times_{M}\mf{X}$.
\end{theorem}
Let $M$ and $\mf{X}$ be as in Theorem \ref{prop:cslrs}.  The proof of Lemma \ref{l:diffqs} goes 
through word for word (after replacing ``an \'etale morphism 
$U\rightarrow\mf{X}$'' by ``a smooth morphism $U\rightarrow\mf{X}$'') to show
\[
j_*\DRB=\epsilon_*f_*\XSB,
\]
where $\epsilon:M_{lis-et}\rightarrow M_{et}$ is the canonical morphism of topoi.  
\begin{remark}
\emph{Since $\epsilon_*f_*\XS^s=\epsilon_*f_*\HH^0(\pLXS)$, the above equality shows that $j_*\DR^s$ is coherent, which is 
not \emph{a\ priori} obvious.}
\end{remark}
To simplify notation, throughout the rest of this subsection we suppress $\epsilon$.  
%As usual, given any stack $\mathcal{Z}$ 
%over $S$, we let $\mathcal{Z}'$ denote the pullback of $\mathcal{Z}$ by $F_S$.  Given any morphism 
%$g:\mathcal{Z}_1\rightarrow\mathcal{Z}_2$ of $S$-schemes, we let $g':\mathcal{Z}'_1\rightarrow\mathcal{Z}'_2$ 
%be the induced morphism.  Relative Frobenii are denoted by $F$.

\begin{proposition}[Cartier isomorphism]
\label{prop:cartieriso}
Let $\mf{X}$ be a smooth tame stack over $S$ which lifts mod $p^2$, and let $f:\mf{X}\rightarrow M$ be its coarse 
space.  Then there is a canonical isomorphism 
\[ 
\HH^t(f'_*F_*\XSB)\stackrel{\simeq}{\rightarrow}f'_*\XPS^t.
\]
If we further assume that $\mf{X}$ and $M$ are as in Theorem \ref{prop:cslrs}, then 
\[
\HH^t(F_*j_*\DRB)\stackrel{\simeq}{\rightarrow}j'_*\DR^t.
\]
\end{proposition}
\begin{proof}
For any left exact functor $G:\mathcal{A}\rightarrow\mathcal{B}$ of abelian categories and any complex $A^\bullet$ of 
objects of $\mathcal{A}$, there is a canonical morphism 
$\HH^t(GA^\bullet)\rightarrow G\HH^t(A^\bullet)$: the map 
$GA^\bullet\rightarrow RGA^\bullet$ induces a morphism from $\HH^t(GA^\bullet)$ to the $E_2^{0t}$-term of 
the spectral sequence $E_2^{st}=R^sG\HH^t(A^\bullet)\Rightarrow R^nG(A^\bullet)$.\\
\\
  For us this yields the (global) map
  $\phi:\HH^t(f'_*F_*\XSB)\rightarrow f'_*\HH^t(F_*\XSB)$ $=f'_*\XPS^t$.
  To prove this is an isomorphism, we need only do so locally.  So, by
  %Lemma~\ref{l:semidirect} 
  \cite[Prop~5.2]{cstlr} and Proposition~\ref{l:Artinfroblift}, we
  are reduced to the case $\mf{X}=[U/G]$, where $U$ is smooth and
  affine, $G=\mathbb{G}_{m,S}^r\rtimes H$ for some finite \'etale
  constant group scheme $H$, and both $\mf{X}$ and the relative
  Frobenius lift mod $p^2$.  Let $X_\bullet$ be the simplical scheme
  obtained by taking fiber products of $U$ over $\mf{X}$, and let
  $a:X_\bullet\rightarrow\mf{X}$ be the augmentation map.  Since
  $U\rightarrow\mf{X}$ lifts mod $p^2$, Theorem~\ref{thm:disimp}
  yields a quasi-isomorphism
\[
\varphi:\bigoplus_{t<p}\Omega^t_{X'_\bullet/S}[-t]\stackrel{\simeq}{\rightarrow}\tau_{<p}F_*\Omega^\bullet_{X_\bullet/S}.
\]
In this local setting, $\varphi$ is a $\emph{morphism\ of\ complexes}$ by 
Remark \ref{rmk:morphcomplex}.  We can therefore apply $(f'a)_*$.  Subsequently taking cohomology, we have a morphism 
$f'_*\XPS^t\stackrel{f'_*\varphi^t}{\longrightarrow} \HH^t(f'_*F_*\XSB)$.  We show that 
\[
\psi:f'_*\HH^t(F_*\XSB)\stackrel{(f'_*C^{-1})^{-1}}{\longrightarrow} f'_*\XPS^t\stackrel{f'_*\varphi^t}{\longrightarrow} 
\HH^t(f'_*F_*\XSB)
\]
and $\phi$ are inverses.  Note that in this local setting $f'_*$ is simply ``take $G$-invariants'', and that 
$\phi:\HH^t((F_*\Omega_U^\bullet)^G)\rightarrow \HH^t(F_*\Omega_U^\bullet)^G$ is $[\alpha]\mapsto (\alpha)$, 
where we use square, resp. round brackets to denote classes in $\HH^t((F_*\Omega_U^\bullet)^G)$, resp. 
$\HH^t(F_*\Omega_U^\bullet)^G$.\\
\\
In general, one does not expect the map $(\alpha)\mapsto [\alpha]$ to be well-defined, but we show here that 
this is precisely what $\psi$ is.  Let $(\omega)\in\HH^t(F_*\Omega_U^\bullet)^G$.  Via the Cartier isomorphism 
$\Omega_{U'}^t\stackrel{(C^{-1})^G}{\longrightarrow} \HH^t(F_*\Omega_U^\bullet)^G$, we know that $(\omega)$ is of 
the form 
\[
(\sum f_{i_1,\dots,i_t} x_{i_1}^{p-1}\dots x_{i_t}^{p-1} dx_{i_1}\wedge\dots\wedge dx_{i_t}),
\]
where
\[
\sum f_{i_1,\dots,i_t} (dx_{i_1}\otimes1)\wedge\dots\wedge (dx_{i_t}\otimes1)\in (\Omega_{U'}^t)^G.
\]
The Deligne-Illusie map $\varphi^q$ sends this $G$-invariant form to
\[
\eta=\sum f_{i_1,\dots,i_t} (x_{i_1}^{p-1}dx_{i_1}+du(x_{i_1}))\wedge\dots\wedge (x_{i_t}^{p-1}dx_{i_t}+du(x_{i_t})),
\]
where $u(x)$ is the reduction mod $p$ of any $y$ satisfying $\til{F}^*(d\til{x}\otimes 1)=\til{x}^p d\til{x} + py$.  
So, $\psi$ sends $(\omega)$ to $(\eta)$.  But since $(du(x))=0$, we see that $\psi$ is the map sending 
$(\alpha)$ to $[\alpha]$.
\end{proof}

\subsection{Degeneracy of Various Spectral Sequences and a Vanishing Theorem}
Let $\mf{X}$ and $M$ be as in Theorem \ref{prop:cslrs}.  Our immediate goal is to show the degeneracy of the 
hypercohomology spectral sequence for $j_*\DRB$ when $\mf{X}$ is proper and lifts mod $p^2$.  If $\XS^1$ were coherent, 
then since $\mf{X}$ is tame, we would have $j_*\DRB=\epsilon_*f_*\XSB=\epsilon_*Rf_*\XSB$.  
The proof of Theorem \ref{thm:st} would then apply directly to show the degeneracy of 
$E_1^{st}=H^t(j_*\DR^s)\Rightarrow H^n(j_*\DRB)$.  Since $\XS^1$ is not coherent, we must take 
a different approach.  As we explain below, the Cartier isomorphism for $j_*\DRB$ proved in the last subsection implies 
that the degeneracy of the above hypercohomology spectral sequence is equivalent to the 
degeneracy of the conjugate spectral sequence $E_2^{st}=H^s(\HH^t(j_*\DRB))\Rightarrow H^n(j_*\DRB)$.  We show that 
this latter spectral sequence degenerates by comparing it to the spectral sequence 
$E_2^{st}=H^s(R^tf_*\XSB)\Rightarrow H^n(\XSB)$ over which we have more control due to the Deligne-Illusie result of 
Section $\ref{sec:simplicial}$.\\
\\
As in the last subsection, we suppress $\epsilon:M_{et}\rightarrow M_{lis-et}$.  The following is the key technical lemma 
we use to prove the degeneracy of the hypercohomology spectral sequence for $j_*\DRB$.
\begin{lemma}
\label{l:L-shaped}
Let $E$ and ${}^{\backprime} E$ be two first quadrant $E_2$ spectral sequences.  Suppose that for $s\neq0$, every differential 
${}^{\backprime} E_r^{st}\rightarrow {}^{\backprime} E_r^{s+r,t-(r-1)}$ is zero.  Suppose further that we are given a morphism $E\rightarrow {}^{\backprime} E$ of 
spectral sequences such that the induced morphism $E_r^{st}\rightarrow {}^{\backprime} E_r^{s+r,t-(r-1)}$ is zero for all $r$, $s$, and $t$, 
and such that $E_2^{st}\rightarrow {}^{\backprime} E_2^{st}$ is an injection for all $s$ and $t$.  Then $E$ degenerates.
\end{lemma}
\begin{proof}
We claim that the morphism $E_r^{st}\rightarrow {}^{\backprime} E_r^{st}$ is an injection for $s\geq r$.  Note that this is 
enough to prove the lemma since for all $s$, the square
\[
\xymatrix{
E_r^{st}\ar[r]\ar[d]_{d_r^{st}} & {}^{\backprime} E_r^{st}\ar[d]\\
E_r^{s+r,t-(r-1)}\ar[r] & {}^{\backprime} E_r^{s+r,t-(r-1)}
}
\]
commutes, the composite is zero, and $E_r^{s+r,t-(r-1)}\rightarrow {}^{\backprime} E_r^{s+r,t-(r-1)}$ is an injection; this shows that all 
differentials $d_r^{st}$ are zero.\\
\\
We now prove the claim by induction.  It is true for $r=2$, so we may assume $r>2$.  Let $s\geq r$ and consider the 
commutative diagram
\[
\xymatrix{
E_{r-1}^{s-(r-1),t+(r-2)}\ar[r]\ar[d] & E_{r-1}^{st}\ar[r]\ar[d] & E_{r-1}^{s+r-1,t-(r-2)}\ar[d]\\
{}^{\backprime} E_{r-1}^{s-(r-1),t+(r-2)}\ar[r] & {}^{\backprime} E_{r-1}^{st}\ar[r] & {}^{\backprime} E_{r-1}^{s+r-1,t-(r-2)}
}
\]
By the inductive hypothesis, all vertical arrows are injective and all arrows on the bottom row are zero.  It follows that 
all arrows on the top arrow are zero, and so $E_r^{st}\rightarrow {}^{\backprime} E_r^{st}$ is injective.
\end{proof}
\begin{theorem}
\label{thm:conj}
Let $\mf{X}$ and $M$ be as in Theorem \ref{prop:cslrs}.  If $M$ has isolated singularities, and $\mf{X}$ is proper and 
lifts mod $p^2$, then the conjugate spectral sequence 
\[
E_2^{st}=H^s(\HH^t(j_*\DRB))\Rightarrow H^n(j_*\DRB)
\]
degenerates for $s+t<p$.
\end{theorem}
\begin{proof}
  %Let $X_\bullet$ be as in Remark~\ref{rmk:liftablecover} and let
  Choose a smooth cover 
  $X\rightarrow\mf{X}$ by a smooth scheme such that the cover lifts
  mod $p^2$, and let $X_\bullet$ be the simplicial scheme obtained 
  by taking fiber products of $X$ with itself over $\mf{X}$.  Let
  $a:X_\bullet\rightarrow\mf{X}$ be the augmentation map.  By Theorem
  \ref{thm:disimp}, we have an isomorphism
  $\bigoplus_{i<p}\Omega^i_{X'_\bullet/S}[-i]\stackrel{\simeq}{\rightarrow}
  \tau_{<p}F_*\Omega^\bullet_{X_\bullet/S}$ in the derived category,
  and therefore, also an isomorphism
\[
\bigoplus_{i<p}R(f'_*a_*)\Omega^i_{X'_\bullet/S}[-i]
\stackrel{\simeq}{\longrightarrow}
\tau_{<p}R(f'_*a_*)F_*\Omega^\bullet_{X_\bullet/S}.
\]
The first of these isomorphisms implies that the Leray spectral
sequence
\[
{}^{\backprime\backprime} E_2^{st}=R^sf'_*\HH^t(F_*\XSB)
\Rightarrow R^nf'_*F_*\XSB
\]
degenerates and that the extension problem is trivial.  The second of
the two isomorphisms shows that the spectral sequence
\[
{}^{\backprime} E_2^{st}=H^s(R^tf'_*F_*\XSB)\Rightarrow H^n(\XSB)
\]
decomposes as the direct sum $\bigoplus {}^{\backprime}_iE$ of Leray
spectral sequences, where
\[
{}^{\backprime}_iE_2^{st}=H^s(R^{t-i}f'_*\XPS^i)\Rightarrow H^{n-i}(\XPS^i).
\]
Note that the morphism $f'_*F_*\XSB\rightarrow Rf'_*F_*\XSB$ induces a
morphism of spectral sequences $E\rightarrow {}^{\backprime} E$, where
\[
E_2^{st}=H^s(\HH^t(f'_*F_*\XSB))\Rightarrow H^n(f_*\XSB).
\]
By the degeneracy of ${}^{\backprime\backprime} E$, the morphism
$\varphi:\HH^t(f'_*F_*\XSB)\rightarrow R^tf'_*F_*\XSB$ factors as
\[
\HH^t(f'_*F_*\XSB)\rightarrow f'_*F_*\XS^t \hookrightarrow R^tf'_*F_*\XSB.
\]
But this first morphism is precisely how the Cartier isomorphism of
Proposition~\ref{prop:cartieriso} was defined.  From this and the
fact that the extension problem for ${}^{\backprime\backprime} E$ is
trivial, we have a split short exact sequence
\[
0\longrightarrow \HH^t(f'_*F_*\XSB)
\stackrel{\varphi}{\longrightarrow} R^tf'_*F_*\XSB\longrightarrow 
\bigoplus_{\substack{i+j=t\\ j>0}}R^jf'_*\XPS^i \longrightarrow 0.
\]
It follows that $E_2^{st}$ is mapped isomorphically to the direct
summand ${}^{\backprime}_tE_2^{st}$ of ${}^{\backprime} E_2^{st}$.
This implies that for all $r$, $s$, and $t$, the induced morphism
$E_r^{st}\rightarrow {}^{\backprime} E_r^{s+r,t-(r-1)}$ is zero.

Note that
\[
j'^*R^tf'_*\XPS^i=j'^*R^tf'_*\bigwedge^i L_{\mf{X}'/S}
=(f^0)'_*\HH^t(\bigwedge^i L_{(\mf{X}^0)'/S})=0.
\]
It follows that $R^tf'_*\XPS^i$ is supported at the singular locus of
$M'$, and since $M$ is assumed to have isolated singularities,
$H^s(R^tf'_*\XS^s)=0$ for $s$ and $t$ positive.  We see then that
${}^{\backprime}_iE_2^{st}$ is zero if $t>i$ and $s>0$, or if $t<i$.
Therefore, the differential ${}^{\backprime} E_r^{st}\rightarrow
{}^{\backprime} E_r^{s+r,t-(r-1)}$ is zero if $s\neq0$.  
>From Lemma~\ref{l:L-shaped}, it follows that $E$ degenerates.
\end{proof}
\begin{remark}
\label{rmk:conjlocfree}
Let $\mathcal{E}$ be a locally free sheaf on $M'$.  Tensoring the isomorphism 
\[
\bigoplus_{i<p}R(f'_*a_*)\Omega^i_{X'_\bullet/S}[-i]\stackrel{\simeq}{\longrightarrow}
\tau_{<p}R(f'_*a_*)F_*\Omega^\bullet_{X_\bullet/S}
\]
with $\mathcal{E}$, we see that the Leray spectral sequence
\[
{}^{\backprime} E_2^{st}=H^s(R^tf'_*F_*\XSB\otimes \mathcal{E})\Rightarrow 
H^n(Rf'_*F_*\XSB\otimes \mathcal{E})
\]
decomposes as the direct sum of spectral sequences.  
%\[
%{}^{\backprime}_iE_2^{st}=H^s(R^{t-i}f'_*\XPS^i\otimes \mathcal{E})\Rightarrow H^{n-i}(Rf'_*\XPS^i\otimes \mathcal{E}).
%\]
The proof of Theorem \ref{thm:conj} then shows that the spectral sequence
\[
E_2^{st}=H^s(\HH^t(f'_*F_*\XSB)\otimes \mathcal{E})\Rightarrow H^n(Rf'_*F_*\XSB\otimes \mathcal{E}).
\]
degenerates for $s+t<p$.
\end{remark}
%Let $X_\bullet$ be as in Remark \ref{rmk:liftablecover} and let $a:X_\bullet\rightarrow\mf{X}$ be the augmentation map.  As 
%in the proof of Theorem \ref{thm:conj}, there is an isomorphism
%\[
%\bigoplus_{i<p}R(f'_*a_*)\Omega^i_{X'_\bullet/S}[-i]\stackrel{\simeq}{\longrightarrow}
%\tau_{<p}R(f'_*a_*)F_*\Omega^\bullet_{X_\bullet/S}
%\]
%in the derived category.  Since $\mathcal{M}^p$ is locally free, tensoring the above isomorphism with 
%$\mathcal{M}^p$ shows that the Leray spectral sequence
%\[
%{}^{\backprime} E_2^{st}=H^s(R^tf'_*F_*\XSB\otimes \mathcal{M}^p)\Rightarrow 
%H^n(Rf'_*F_*\XSB\otimes \mathcal{M}^p)
%\]
%decomposes as the direct sum $\bigoplus {}^{\backprime}_iE$ of Leray spectral sequences, where
%\[
%{}^{\backprime}_iE_2^{st}=H^s(R^{t-i}f'_*\XPS^i)\Rightarrow H^{n-i}(Rf'_*\XPS^i\otimes \mathcal{M}^p).
%\]
%The morphism $f'_*F_*\XSB\rightarrow Rf'_*F_*\XSB$ induces a morphism of spectral sequences 
%$E\rightarrow {}^{\backprime}E$, where 
%\[
%E_2^{st}=H^s(\HH^t(f'_*F_*\XSB)\otimes \mathcal{M}^p)\Rightarrow H^n(Rf'_*F_*\XSB\otimes \mathcal{M}^p).
%\]
%The Cartier isomorphism of Proposition \cite{prop:cartieriso} then shows that $E_2^{st}$ is mapped isomorphically to 
%${\backprime}_tE_2^{st}$.
\begin{theorem}
\label{thm:hyp}
Let $\mf{X}$ and $M$ be as in Theorem \ref{prop:cslrs}.  If $M$ has isolated singularities, and $\mf{X}$ is proper and 
lifts mod $p^2$, then the hypercohomology spectral sequence
\[
E_1^{st}=H^t(j_*\DR^s)\Rightarrow H^n(j_*\DRB)
\]
degenerates for $s+t<p$.
\end{theorem}
\begin{proof}
By the Cartier isomorphism 
\[
H^s(\HH^t(j_*\DRB))=H^s(\HH^t(f'_*F_*\XSB))=H^s(f'_*\XPS^t).
\]
But $H^s(f'_*\XPS^t)=H^s(f_*\XS^t)\otimes_{k,F_k}k$; in particular, 
\[
\dim_k H^s(\HH^t(j_*\DRB))=\dim_k H^s(f_*\XS^t).
\]
By Corollary \ref{cor:tamedegen}, the above cohomology groups are finite-dimensional 
$k$-vector spaces.  The degeneracy of the conjugate spectral sequence shows
\[
H^n(j_*\DRB)\simeq \bigoplus_{s+t=n}H^s(\HH^t(j_*\DRB)),
\] 
and so
\[
\dim_k H^n(j_*\DRB)=\sum_{s+t=n}\dim_k H^s(f_*\XS^t),
\]
which implies the degeneracy of the hypercohomology spectral sequence by Lemma \ref{l:homological}.
\end{proof}
Although our proof of Theorem $\ref{thm:st}$ goes through stacks, the statement of the theorem is purely 
scheme-theoretic.  We would similarly like to remove the stack from the statement of Theorem $\ref{thm:hyp}$.  We 
can do so when $M$ has large enough dimension.
\begin{theorem}
\label{thm:lift}
Let $M$ be a proper $S$-scheme with isolated linearly reductive singularities.  If $\dim M\geq 4$ and $M$ lifts 
mod $p^2$, then the stack $\mf{X}$ of Theorem \ref{prop:cslrs} lifts mod $p^2$, and so
\[
E_1^{st}=H^t(j_*\DR^s)\Rightarrow H^n(j_*\DRB)
\]
degenerates for $s+t<p$.
\end{theorem}
\begin{proof}
Let $m=\dim M$.  
%and let $\mf{X}$ be as in Theorem \ref{prop:cslrs}.  If we can prove $\mf{X}$ lifts mod $p^2$, then we are done.  
The exact triangle 
\[
Lf^*L_{M/S}\longrightarrow L_{\mf{X}/S}\longrightarrow L_{\mf{X}/M}
\]
gives rise to the long exact sequence 
\[
\dots\longrightarrow \Ext^2(L_{\mf{X}/M},\OO_\mf{X})\longrightarrow \Ext^2(L_{\mf{X}/S},\OO_\mf{X})
\longrightarrow \Ext^2(Lf^*L_{M/S},\OO_\mf{X})\longrightarrow\dots. 
%\Ext^3(L_{\mf{X}/M},\OO_\mf{X})\longrightarrow\dots.
\]
Note that
\[
RHom(Lf^*L_{M/S},\OO_\mf{X})=RHom(L_{M/S},Rf_*\OO_\mf{X})=RHom(L_{M/S},\OO_M)
\]
since $Rf_*\OO_\mf{X}=f_*\OO_\mf{X}$ by tameness and $f_*\OO_\mf{X}=\OO_M$ by Keel-Mori \cite{km}.  Since the obstruction to 
lifting $\mf{X}$ lies in $\Ext^2(L_{\mf{X}/S},\OO_\mf{X})$, we need only show $\Ext^2(L_{\mf{X}/M},\OO_\mf{X})=0$.  
We in fact prove $\mathcal{E}xt^s(L_{\mf{X}/S},\OO_\mf{X})=0$ for $s\leq m-2$.\\
\\
Since $(j^0)^*L_{\mf{X}/M}=L_{\mf{X}^0/M^0}=0$, we see
\[
0=Rj^0_*R\mathcal{H}om((j^0)^*L_{\mf{X}/M},\OO_{\mf{X}^0})=R\mathcal{H}om(L_{\mf{X}/M},Rj^0_*\OO_{\mf{X}^0}).
\]
A local cohomology argument given below will show $R^tj^0_*\OO_{\mf{X}^0}\neq 0$ if and only if $t=0,m-1$.  Assuming 
this for the moment, let us complete the proof.  We have a spectral sequence 
\[
E_2^{st}=R^s\mathcal{H}om(L_{\mf{X}/M},R^tj^0_*\OO_{\mf{X}^0})\Rightarrow 
R^n\mathcal{H}om(L_{\mf{X}/M},Rj^0_*\OO_{\mf{X}^0})=0.
\]
The only page with non-zero differentials, then, is the $m^{th}$.  Since $L_{\mf{X}/M}$ is concentrated in 
degrees at most 1, $R^s\mathcal{H}om(L_{\mf{X}/M},R^tj^0_*\OO_{\mf{X}^0})=0$ for $s<-1$.  It follows that 
\[
R^s\mathcal{H}om(L_{\mf{X}/M},j^0_*\OO_{\mf{X}^0})=0
\]
for $s\leq m-2$, which proves the theorem since $j^0_*\OO_{\mf{X}^0}=\OO_{\mf{X}}$.\\
\\
We now turn to the local cohomology argument.  To prove $R^tj^0_*\OO_{\mf{X}^0}\neq 0$ if and only if $t=0,m-1$, we 
can make an \'etale base change.  We can therefore assume $\mf{X}=[U/G]$, where $U$ is smooth and affine, and $G$ is 
%smooth, affine, and 
finite linearly reductive.  Since $M$ has isolated singularities, we can further assume 
$U^0=U\backslash\{x\}$, where $U^0$ is the pullback
\[
\xymatrix{
U^0\ar[r]^i\ar[d]_h & U\ar[d]^g\\
\mf{X}^0\ar[r]^{j^0} & \mf{X}
}
\]
%Since $g$ is a smooth cover, 
The following lemma, then, completes the proof.
\end{proof}
\begin{lemma}
Let $U$ be a normal affine scheme of dimension $m$ and let $x\in U$ be Cohen-Macaulay.  If $U^0=U\backslash\{x\}$ 
and $i:U^0\hookrightarrow U$ is the inclusion, then $R^ti_*\OO_{U^0}\neq 0$ if and only if $t=0,m-1$.
\end{lemma}
\begin{proof}
Note that $R^ti_*\OO_{U^0}$ is the skyscraper sheaf $H^t(\OO_{U^0})$ at $x$.  By normality, $H^0(\OO_{U^0})=H^0(\OO_U)$.  
Since $U$ is affine, the long exact sequence 
\[
\dots \longrightarrow H^n_x(\OO_U)\longrightarrow H^n(\OO_U)\longrightarrow H^n_x(\OO_{U^0})\longrightarrow 
H^n_x(\OO_U)\longrightarrow \dots
\]
shows $H^t(\OO_{U^0})=H^{t+1}_x(\OO_U)$ for $t>0$.  Since $x$ is Cohen-Macaulay, $H^{t+1}_x(\OO_U)\neq0$ if and only 
if $t+1=m$.
\end{proof}

We now prove an analogue of \cite[Lemma 2.9]{di} which Deligne and Illusie use to deduce Kodaira Vanishing.
\begin{lemma}
\label{l:KVlemma}
Let $\mf{X}$ and $M$ be as in Theorem \ref{prop:cslrs}.  Suppose $M$ has isolated singularities, and $\mf{X}$ is proper 
and lifts mod $p^2$.  Let $d$ be the dimension of $M$ and let $N$ be an integer such that $N\leq\inf(d,p)$.  If 
$\mathcal{M}$ is an invertible sheaf on $M$ such that
\[
H^t(j_*\Omega^s_{M^0/S}\otimes\mathcal{M}^p)=0
\]
for all $s+t<N$, then
\[
H^t(j_*\Omega^s_{M^0/S}\otimes\mathcal{M})=0
\]
for all $s+t<N$.
\end{lemma}
\begin{proof}
Let $\mathcal{M}'$ be the pullback of $\mathcal{M}$ to $M'$.  Since $F^*\mathcal{M}'=\mathcal{M}^p$, the projection 
formula shows
\[
H^t(j_*\Omega^s_{M^0/S}\otimes \mathcal{M}^p)=H^t(f'_*F_*\XS^s\otimes \mathcal{M}').
\]
From the hypercohomology spectral sequence
\[
E_2^{st}=H^t(f'_*F_*\XS^s\otimes \mathcal{M}')\Rightarrow H^n(f'_*F_*\XSB\otimes \mathcal{M}'),
\]
we see that $H^n(f'_*F_*\XSB\otimes \mathcal{M}')=0$ for all $n<N$.  Remark \ref{rmk:conjlocfree} shows that the Leray 
spectral sequence
\[
E_2^{st}=H^s(f'_*\XPS^t\otimes\mathcal{M}')\Rightarrow H^n(f'_*F_*\XSB\otimes \mathcal{M}')
\]
degenerates, and so $H^s(f'_*\XPS^t\otimes\mathcal{M}')=0$ for all $s+t<N$.  Since 
\[
\dim_k H^s(f'_*\XPS^t\otimes\mathcal{M}')=\dim_k H^s(j_*\Omega_{M^0/S}^t\otimes\mathcal{M}),
\]
the lemma follows.
\end{proof}

%%%%%%%%%%%%%%
\begin{comment}
\begin{lemma}
\label{l:newcmlemma}
If $M$ is a scheme with linearly reductive singularities, then $j_*\Omega^s_{M^0/S}$ is Cohen-Macaulay for all $s$.
\end{lemma}
\begin{proof}
To prove the statement, it suffices to look \'etale locally, where we can assume $M=U/G$ for a smooth scheme $U$ and finite linearly reductive group scheme $G$.  By Theorem \ref{prop:cslrs}, we can further assume that $U\rightarrow M$ is a $G$-torsor over $M^0$.\\
\\
Since ?????????, $j_*\Omega^s_{M^0/S}=(\Omega^s_{U/S})^G$.  Since $G$ is linearly reductive, the natural inclusion
\[
(\Omega^s_{U/S})^G\longrightarrow \Omega^s_{U/S}
\]
has a splitting, and so $(\Omega^s_{U/S})^G$ is a direct summand of $\Omega^s_{U/S}$.  As $U$ is smooth, $\Omega^s_{U/S}$ is a Cohen-Macaulay $\OO_U$-module, and so $(\Omega^s_{U/S})^G$ is as well.  Since $M$ is Cohen-Macaulay and $U$ is finite over $M$, it follows that $\OO_U$ is a Cohen-Macaulay $\OO_M$-module.  Hence, $(\Omega^s_{U/S})^G$ is a Cohen-Macaulay $\OO_M$-module.
\end{proof}
Using the above two lemmas, we obtain the following vanishing result.
\end{comment}
%%%%%%%%%%%%%%%%%%%%
Unfortunately, we cannot quite deduce from Lemma \ref{l:KVlemma} a general Kodaira Vanishing result.  Following Deligne and Illusie, we would like to show that if $M$ is a projective scheme of dimension $d$ with isolated linearly reductive singularities and $\mathcal{L}$ is an ample line bundle on $M$, then $H^t(j_*\Omega_{M^0/S}^s\otimes\mathcal{L}^{-p^m})=0$ for $m$ sufficiently large.  Lemma \ref{l:KVlemma} would then imply that $m$ can be taken to be 1.  The issue is that the vanishing of these cohomology groups for $m$ large enough is not clear.  Under certain hypothesis, however, we obtain a vanishing theorem.
%This follows from Serre Duality when $j_*\Omega_{M^0/S}^s$ are locally free.  This assumption, however, is too 
%restrictive.  Under the weaker hypothesis that the $j_*\Omega_{M^0/S}^s$ are Cohen-Macaulay, we can still obtain a vanishing 
%result, as is shown in the following proposition.
\begin{proposition}
\label{prop:kv}
Let $M$ be a projective scheme of dimension $d$ with isolated linearly reductive singularities.  Let $\mathcal{L}$ be an ample line bundle on $M$.  If the stack $\mf{X}$ of Theorem \ref{prop:cslrs} lifts mod $p^2$ and 
%, then
the $j_*\Omega_{M^0/S}^s$ are Cohen-Macaulay for all $s$, then
\[
H^t(j_*\Omega_{M^0/S}^s\otimes\mathcal{L}^{-1})=0
\]
for all $s+t<\inf(d,p)$.
\end{proposition}
\begin{proof}
By Lemma \ref{l:KVlemma}, we need only prove that $H^t(j_*\Omega_{M^0/S}^s\otimes\mathcal{L}^{-p^m})=0$ for $m$ 
sufficiently large.  Grothendieck Duality shows 
\[
H^t(j_*\Omega_{M^0/S}^s\otimes\mathcal{L}^{-p^m})^\vee=
\Ext^{d-t}(j_*\Omega_{M^0/S}^s\otimes\mathcal{L}^{-p^m},\omega_M).
\]
Since the $\mathcal{E}xt^{d-t}(j_*\Omega_{M^0/S}^s,\omega_M)$ are coherent, the local-global Ext spectral sequence 
shows that for $m$ sufficiently large, 
\[
H^t(j_*\Omega_{M^0/S}^s\otimes\mathcal{L}^{-p^m})^\vee=
\Gamma(\mathcal{E}xt^{d-t}(j_*\Omega_{M^0/S}^s,\omega_M)\otimes\mathcal{L}^{p^m}).
\]
For all $x\in M$, 
\[
\mathcal{E}xt^{d-t}(j_*\Omega_{M^0/S}^s,\omega_M)_x=\Ext^{d-t}_{\OO_x}((j_*\Omega_{M^0/S}^s)_x,\omega_{\OO_x}).
\]
Since $M$ and the $j_*\Omega_{M^0/S}^s$ are Cohen-Macaulay, 
% by Lemma \ref{l:newcmlemma}, and since $M$ is Cohen-Macaulay, 
local duality shows that for $t<d$,
\[
\mathcal{E}xt^{d-t}(j_*\Omega_{M^0/S}^s,\omega_M)=0,
\]
thereby completing the proof.
\end{proof}

We conclude by showing that the hypercohomology spectral sequence 
\[
E_1^{st}=R^tf'_*F_*\XS^s\Rightarrow R^nf'_*F_*\XSB
\]
degenerates at $E_2$ and that the only potentially non-zero differentials on the first page are those on the zero-th row.
\begin{lemma}
\label{l:general}
Let $\mathcal{A}$ and $\mathcal{B}$ be abelian categories and let $F:\mathcal{A}\rightarrow\mathcal{B}$ be a left exact 
functor.  Suppose that $\mathcal{A}$ has enough injectives.  If $A^\bullet$ is a complex of objects in $\mathcal{A}$ and 
$C^\bullet$ denotes the cone in the derived category $D(\mathcal{A})$ of the canonical morphism 
$FA^\bullet\rightarrow RFA^\bullet$, then there is a spectral sequence 
\[
{}^{\backprime} E_1^{st}=\left\{
\begin{array}{lr}
R^tFA^s & t>0\\
0 & t=0
\end{array}
\right.
\Rightarrow \HH^n(C^\bullet).
\]
If in the hypercohomology spectral sequence $E_1^{pq}=R^qFA^p\Rightarrow R^nFA^\bullet$, the differentials 
$E_r^{s-r,r-1}\rightarrow E_r^{s,0}$ are zero for all $r\geq2$, then for every $n$,
\[
0\rightarrow \HH^n(FA^\bullet)\rightarrow R^nFA^\bullet\rightarrow \HH^n(C^\bullet)\rightarrow 0
\]
is a short exact sequence.
\end{lemma}
\begin{proof}
The existence of the spectral sequence ${}^{\backprime} E$ is shown as follows.  Let $A^s\rightarrow I^{s,\bullet}$ be an 
injective resolution of $A^s$.  The cone $C^\bullet$ is then quasi-isomorphic to the total complex of 
\[
%\textrm{Tot}\left(
\begin{array}{ccc}
\vdots & \vdots & \\
FI^{01} & FI^{11} & \dots\\
FI^{00} & FI^{10} & \dots\\
FA^0 & FA^1 & \dots
\end{array}
%\right)
\]
where $FA^0$ has bidegree $(-1,0)$.  The spectral sequence associated to this double complex in which we begin by taking 
cohomology vertically is our desired ${}^{\backprime}E$.\\
\\
Note that there is a morphism of spectral sequences $E\rightarrow {}^{\backprime}E$.  If the differentials 
$E_r^{s-r,r-1}\rightarrow E_r^{s,0}$ are zero for all $r\geq2$, then the morphism of spectral sequences induces an 
isomorphism $E_\infty^{st}\stackrel{\simeq}{\rightarrow} {}^{\backprime}E_\infty^{st}$ for $t\neq0$.  It follows that 
$\HH^n(C^\bullet)$ is equal to $R^nFA^\bullet$ modulo the bottom part of its filtration, namely 
$E_\infty^{n0}=\HH^n(FA^\bullet)$.
\end{proof}
\begin{proposition}
Let $\mf{X}$ and $M$ be as in Theorem \ref{prop:cslrs}.  If $M$ has isolated singularities, and $\mf{X}$ is proper and 
lifts mod $p^2$, then the hypercohomology spectral sequence 
\[
E_1^{st}=R^tf'_*F_*\XS^s\Rightarrow R^nf'_*F_*\XSB
\]
degenerates at $E_2$, and for $t\neq0$, the differentials $E_1^{st}\rightarrow E_1^{s+1,t}$ are zero.
\end{proposition}
\begin{proof}
Let $C^\bullet$ be the cone of the canonical morphism $f'_*F_*\XSB\rightarrow Rf'_*F_*\XSB$.  Note that for $t>0$, we have
\[
j'^*R^tf'_*F_*\XS^s=F_*j^*R^tf_*\pLXS=(f^0)'_*\HH^t(\bigwedge^s L_{(\mf{X}^0)'/S})=0,
\]
and so $R^tf'_*F_*\XS^s$ is supported at the singular 
locus of $M'$; in particular, the $R^tf'_*F_*\XS^s$ are torsion.  On the other hand, 
$R^0f'_*F_*\XS^s=F_*j_*\Omega^s_{M^0/S}$, which is reflexive, and hence torsion-free.  As a result, for $r\geq2$ every 
differential $E_r^{s-r,r-1}\rightarrow E_r^{s,0}$ is zero, and $E_r^{st}$ is supported at the singular locus of $M'$ for all 
$t\neq0$ and all $s$ and $r$.  So, to prove the proposition, we need only show that the spectral sequence 
\[
{}^{\backprime} E_1^{st}=\left\{
\begin{array}{lr}
R^tf'_*F_*\XS^s & t>0\\
0 & t=0
\end{array}
\right.
\Rightarrow \HH^n(C^\bullet)
\]
of Lemma \ref{l:general} degenerates.\\
\\
Since $M$ is assumed to have isolated singularities, for any short exact sequence
\[
0\rightarrow \mathcal{F} \rightarrow \mathcal{G} \rightarrow \mathcal{Q} \rightarrow 0
\]
with $\mathcal{F}$ supported at the singular locus,
\[
0\rightarrow \Gamma(\mathcal{F}) \rightarrow \Gamma(\mathcal{G}) \rightarrow \Gamma(\mathcal{Q}) \rightarrow 0
\]
is short exact as well.  Furthermore, $\Gamma(\mathcal{F})=\bigoplus_{x\in M}\mathcal{F}_x$, so $\mathcal{F}$ is zero if and 
only if $\Gamma(\mathcal{F})$ is zero.  It follows that we have a spectral sequence 
\[
{}^{\backprime\backprime}E_1^{st}=\left\{
\begin{array}{lr}
\Gamma(R^tf'_*F_*\XS^s) & t>0\\
0 & t=0
\end{array}
\right.
\Rightarrow \Gamma(\HH^n(C^\bullet))
\]
whose degeneracy is equivalent to that of ${}^{\backprime}E$.  By Lemma \ref{l:general}, there is a short exact sequence 
\[
0\longrightarrow \HH^n(f'_*F_*\XSB)\longrightarrow R^nf'_*F_*\XSB\longrightarrow \HH^n(C^\bullet)\longrightarrow 0.
\]
Comparing this with the short exact sequence
\[
0\longrightarrow \HH^t(f'_*F_*\XSB)\stackrel{\varphi}{\longrightarrow} R^tf'_*F_*\XSB\longrightarrow 
\bigoplus_{\substack{i+j=t\\ j>0}}R^jf'_*\XPS^i \longrightarrow 0
\]
proved in Theorem \ref{thm:conj}, we see 
\[
\HH^n(C^\bullet)=\bigoplus_{\substack{s+t=n\\ t>0}}R^tf'_*\XPS^s.
\]
It follows that
\[
\sum_{s+t=n}\dim_k {}^{\backprime\backprime} E_1^{st}=\dim_k\Gamma(\HH^n(C^\bullet))
\]
which shows the degeneracy of ${}^{\backprime\backprime}E$ by Lemma \ref{l:homological}.
\end{proof}

\end{document}